\documentclass [12pt,oneside] {amsart}
\usepackage{latexsym,amsfonts,amssymb,euler,euscript,amscd,fullpage}
\usepackage[pdftex]{graphicx}

\theoremstyle{plain}
\newtheorem {Thm} {Theorem}[section]
\newtheorem* {Thm*} {Theorem}
\newtheorem* {Prop*} {Proposition}
\newtheorem {Lem}[Thm] {Lemma}
\newtheorem {Prop}[Thm] {Proposition}
\newtheorem {Cor}[Thm] {Corollary}

\newtheorem {Conj}[Thm] {Conjecture}
\theoremstyle {definition}
\newtheorem {Def}[Thm] {Definition}
\newtheorem {Cons}[Thm] {Construction}
\newtheorem {Que}[Thm] {Question}

\newtheorem {Rem}[Thm] {Remark}
\newtheorem {Exa}[Thm] {Example}
\newenvironment{Pf}[1]{{\noindent\sc Proof #1:}}{\qed\\}
\newcommand {\specialmap} [4] {\text {$ #1\negmedspace : #2 #3 #4 $}}
\newcommand {\map} [3] {\specialmap {#1} {#2}{\to} {#3}}

\newcommand {\isomap} [3] {\specialmap {#1} {#2}{\overset {\cong{\phantom{.}}} 
            {\longrightarrow}} {#3}}

\newcommand {\Aut}{\operatorname{Aut}}
\newcommand {\Hom} {\operatorname {Hom}}
\newcommand {\id} {\operatorname{id}}

\newcommand {\at}[1] {\arrowvert_{#1}}
\newcommand {\del} {\partial}
\newcommand {\ind} {\operatorname{ind}}
\renewcommand {\(} {\left(}
\renewcommand {\)} {\right)} 

\renewcommand {\geq} {\geqslant}

\newcommand {\sub} {\subseteq}

\newcommand {\CC} {\mathbb C}
\def\Pfeil{{\unitlength 1em\begin{picture}(0,.1)
\put(0,.1){\vector(1,0){1.5}}\end{picture}}}
\def\colimname{{\unitlength.1em
\raisebox{-2.7\unitlength}{\begin{picture}(15.5,9.5)(0,0)
\put(0,2.7){$\operatorname{lim}$}
\put(.05,-.1){\Pfeil} 
\end{picture}}}}
\def\Colim{\mathop{\colimname}}
\newcommand {\on}[1] {\operatorname{#1}}

\newcommand {\BG} {{\on{BG}}}

\newcommand {\BZZ} {{\on{B\ZZ^2}}}

\newcommand {\dual} {{\check{\phantom x}}}

\renewcommand {\gg} {{(g_1,\dots,g_h)}}
\newcommand {\GG}[1] {{\frac{1}{|#1|}}}

\newcommand {\HH} {{\mathfrak H}}

\newcommand {\HomAE} {{\Hom(\hat A,E)}}
\renewcommand {\ind}[2] {\on{ind}\arrowvert_{#1}^{#2}} 
\newcommand {\Inj} {\operatorname{Inj}}

\newcommand {\HKR} {Hopkins-Kuhn-Ravenel }

\renewcommand {\leq} {\leqslant}

\newcommand {\Lt} {{\langle\tau,1\rangle}}
\newcommand {\Ltp} {{\langle\tau',1\rangle}}

\newcommand {\mmod} {/\!\!/}

\newcommand {\ob} {{\operatorname{ob}}}

\newcommand {\orb} {{\on{orb}}}

\newcommand {\pt} {\on{pt}}

\newcommand {\res}[2] {\on{res}\arrowvert^{#1}_{#2}}

\newcommand {\SL} {{\on{SL}_2(\ZZ)}}
\newcommand {\Sn} {{\Sigma_n}}
\newcommand {\Sm} {{\Sigma_m}}
\newcommand {\Snm} {{\Sigma_{n+m}}}

\newcommand {\Stab} {\on{Stab}}

\newcommand {\Sym} {\on{Sym}}
\newcommand {\tensor}{\otimes}

\newcommand {\tHn} {\on{\widehat T_n}}

\newcommand {\TT} {{\T\in\mathcal{T}}}

\newcommand {\inj} {{\on{inj}}}

\newcommand {\XinjB} {{X_B^\inj}}

\newcommand {\ZZn} {{\ZZ/n\ZZ}}
\newcommand {\ZZl} {{\ZZ/l\ZZ}}

\newcommand {\ZZ} {\mathbb Z}

\newcommand {\lps} {[\! [}
\newcommand {\rps} {]\! ]}

\newcommand {\ps}[1] {\lps #1\rps}

\usepackage [all]{xy}
\SelectTips{cm}{12}
\title{Hecke operators in equivariant elliptic cohomology and
  generalized Moonshine} 
\author{Nora Ganter\\ University of Illinois at Urbana Champaign\\ \\ 
  dedicated to John McKay}
\thanks{The author was supported by NSF-grant DMS-0504539. During her
  work on this paper, the author enjoyed the hospitality of Concordia
  University and of the Mathematical Sciences Research Institute}
\date {\today}
\renewcommand{\check} {^\vee}
\newcommand {\cyl}[2] {{\on{cyl}_{#1}(#2)}}
\newcommand {\Dev}{{\on{Dev}}}
\newcommand {\Ell}{{\on{Ell}}}

\renewcommand {\hat}{\widehat}
\newcommand {\inv}{^{-1}}
\newcommand {\mA} {{\mathcal A}}
\newcommand {\mC} {{\mathcal C}}
\newcommand {\mF} {{\mathcal F}}
\newcommand {\mH} {{\mathcal H}}
\newcommand {\mL} {{\mathcal L}}
\newcommand {\mM} {{\mathcal M}}
\newcommand {\MconnA}{{\mM^{\on{conn}}_A}}
\newcommand {\MzeroA}{{\widetilde{\mM}_{0}(A)}}

\newcommand {\ol}[1] {\overline{#1}}
\newcommand{\PG}{{\overline{\mathcal C}(G)}}
\newcommand{\Prin}[2]{{\on{Prin}_{#1}(#2)}}
\renewcommand {\gg} {{\langle g\rangle}}
\newcommand {\sgn} {{\on{sgn}}}
\newcommand {\sgnst} {{\sgn(\sigma,\rho)}}
\newcommand {\st} {{(\sigma,\rho)}}
\newcommand {\bbS}{{\mathbb S}}
\renewcommand {\Sym}[1]{{\on{Sym}^{(2)}_{#1}}}
\newcommand {\sym}[1]{{\on{sym}^{(2)}_{#1}}}
\newcommand {\tr}[2] {{\on{Tr}(#1 | #2)}}
\newcommand {\trs}[1]{{\tr{\sigma}{#1}}}
\renewcommand {\TT}{{\mathbb T}}

\newcommand{\Witten} {{W}}
\newcommand{\ZxZ} {{\ZZ^2}}
\begin{document}
%
\maketitle
\tableofcontents
\section{Introduction}
\subsection{Background on generalized Moonshine}
We recall Norton's generalized Moonshine conjecture.
Let $g,h$ be a pair of commuting elements in the Fischer-Griess
Monster group $M$. We denote the centralizer of $g$ in $M$ by $C_M(g)$.
\begin{Conj}[{Norton \cite{Norton}}]
  It is possible to associate with each such pair a modular function
  $f(g,h;\tau)$ with the following properties:
  \begin{enumerate}
    \renewcommand{\labelenumi}{(\alph{enumi})}
    \item Up to a constant factor (which will be a root of
      unity), there is an equality 
      $$f(g^ah^c,g^bh^d;\tau) = f\(g,h,\frac{a\tau+b}{c\tau+d}\)$$ whenever
      $ad-bc = 1$, 
    \item For any group element $g$ and nonzero rational number $l$
      the coefficient of $q^l = e^{2\pi i l \tau}$ in $f(g,h;\tau)$ is, as a
      function of $h$, a character of a central extension\footnote{For
      a discussion of how characters of central extensions of a group
      correspond to twisted characters of this group, see for example
      \cite{Willerton}.} 
      of $C_M(g)$. Note that nonzero
      characters can occur for nonintegral $l$, but that generalized
      characters are not needed.
    \item Conjugation of $g$ and $h$ leaves the function unchanged. 
    \item Unless $f(g,h;\tau)$ is a constant function, its invariance
    group will be a modular group of genus zero, commensurable with
    the standard modular group $\SL$.
  \end{enumerate}
\end{Conj}
Norton further explains: ``It is the simultaneous action of $\SL$
on $\langle g,h\rangle$ and $\tau$ that seems to suggest an explanation
of the Moonshine phenomenon; but this only gives a clue to the actions
of the elements of $\SL$ on $\tau$, so the reason why the full
invariance group should have genus zero remains obscure.''

Such simultaneous actions are well known to topologists: they arise in
the context of equivariant elliptic cohomology. There are more
indications of a connection between Moonshine and elliptic cohomology:

Generalized Moonshine is supposed to specialize to
``classical'' Moonshine for pairs of the form $(1,g)$. 
This part is often referred to as the ``untwisted sector''. 
In this situation, part (b) of Norton's conjecture involves only
integral powers of $q$. Moreover, no central extension is needed, 
i.e., the coefficients are honest characters of $M=C_M(1)$.
%
In this context, the Moonshine functions are usually denoted
$$
  j_{\langle g\rangle}(\tau) = f(1,g;\tau) = q^{-1}+a_1q+a_2q^2+\dots,
$$
since they are generalizations of the $j$-function 
$$
  j_{\langle 1\rangle} = j-744 = q\inv+196884q+21493760 q^2 +\dots
$$
and depend only
on the cyclic subgroup generated by $g$. 

We recall from \cite{Norton} the notions of {\em replicable function}
and {\em twisted Hecke operator}. Let $f$ be a Laurent series in $q$ of
the form
\begin{equation}
  \label{ls-Eqn}
  f(q) = q^{-1} + a_1 q + a_2q^2 + \dots.  
\end{equation}
Then the $n^{th}$ {\em Faber polynomial} $\Phi_n$ of $f$ is the unique
polynomial of degree $n$
for which the Laurent series $\Phi_n(f)$ has the form
$$
  \Phi_n(f) = q^{-n}+b_1q+b_2q^2+\dots.
$$
The series $f$ is called {\em replicable}, if there exist {\em
  replicates} $f^{(a)}$ of $f$, which satisfy the equalities
\begin{equation}
  \label{Faber-Hecke-Eqn}
  \Phi_n(f) = \sum_{\stackrel{ad=n}{0\leq b<d}}f^{(a)}\(\frac{a\tau
    +b}{d}\)
\end{equation}
for all $n\in\mathbb N$.
One shows by induction that the $f^{(a)}$ are uniquely determined if
they exist: assume we know that $f^{(a)}$ is defined for
$a<n$. By solving the equality \eqref{Faber-Hecke-Eqn} for
$f^{(n)}(n\tau)$, we arrive 
at the question whether a series in $q$, which we already defined, is
indeed a series in $q^n$.

A replicable series $f$ is called {\em completely replicable} if all
the replicates of $f$ are again replicable and their replicates are
given by 
$$\(f^{(a)}\)^{(b)}=f^{(ab)}.$$
Complete replicability was proved to be equivalent to the genus
zero property (part (d) of the conjecture) in \cite{Cummins:Gannon}
and \cite{Martin}. 
 
If $f$ is replicable, then the {\em
  $n^{th}$ twisted Hecke operator} acting on $f$ is defined as
\begin{equation}
  \label{twisted-Hecke-Eqn}
  \hat T_n(f) := \frac 1n\Phi_n(f).  
\end{equation}
The classical Moonshine
functions are known to be completely replicable, with the replicates
given by the 
Adams operations
$$
  j^{(a)}_{\langle h\rangle}(q) = j_{\langle h^a\rangle}(q).
$$
The collection of all classical Moonshine functions is a
{\em McKay-Thompson series}, i.e., it is a series of the form (\ref{ls-Eqn}) with
coefficients in the representation ring $R(M)$. 
A function $f$ satisfying condition (b) of Norton's conjecture is
often referred to as a generalized McKay-Thompson series.
The following
definition of replicability for McKay-Thompson series can be found in
\cite{Norton}. It is stricter than asking for the individual functions $j_\gg$
to be replicable but seems to be the correct notion in the context of
Moonshine. 
\begin{Def}\label{untwisted-replicability-Def}
  A McKay-Thompson series $f$ is called replicable, if the formula
 (\ref{Faber-Hecke-Eqn}) holds when $f^{(a)}$ is the $a^{th}$ Adams
 operation applied to $f$.  
\end{Def}
We will rephrase Definition \ref{untwisted-replicability-Def}
in terms of power operations in elliptic cohomology in Section
\ref{replicability-Sec} below. 
\subsection{Background on equivariant elliptic cohomology and
  statement of results}
Adams operations were originally defined in the context of (equivariant)
$K$-theory. Hecke operators in (equivariant) elliptic cohomology are
a generalization\footnote{More precisely, they are the chromatic level
$2$ analogue of the Adams operations (cf.\ also \cite{Ganter:thesis}).}
of the Adams operations \cite{Ando:thesis}, and it
seems a natural question to ask for a ``topological'' explanation for
the relationship between Hecke and Adams operations occurring in Moonshine. We
will start by explaining the geometry behind the generalized Moonshine 
conjecture and its relationship to equivariant elliptic cohomology.
This discussion is not specific to the Monster and can be formulated
for any finite group $G$.
We will translate 
parts (a) and (c) of Norton's conjecture to 

\medskip
``{\em $f$ is an
  element of the twisted equivariant elliptic cohomology of the one point
  space,
  in the sense of 
  Ginzburg, Kapranov and Vasserot: $f\in \Ell_G^\alpha(\pt)$.}''
\medskip

Here $\alpha$ is a cocycle representing an element of 
$$H^3(\BG;\mathbb R/\ZZ) = H^4(\BG, \ZZ).$$ 
In other words, 

\medskip
``{\em $f$ is a section of the line bundle $\mL^\alpha$ over the moduli
  stack of principal $G$-bundles over (complex) elliptic curves.}''
\medskip

Here $\mL^\alpha$ is a line bundle associated to $\alpha$ (cf.\
Section \ref{Freed-Quinn-Sec}).

In this context, Ginzburg, Kapranov and Vasserot
have defined the action of correspondences. 
In particular, one can consider the action of the Hecke correspondence. 
In Section \ref{Hecke-correspondence-Sec}, I will
spell out this geometric definition of the Hecke operators: they act on 
such sections via pull-backs along isogenies.

We return to Norton's conjecture. One easily translates part (b) to 
\medskip

``{\em $f$ is an element of the (twisted)
  Devoto equivariant Tate $K$-theory of the one point space.}'' 

\medskip

Devoto's equivariant Tate
$K$-theory 
$$
  K_{\Dev,g}(X) = \bigoplus_{[g]}K_{C_g}(X^g)\ps{q^\frac1{|g|}}
$$
was introduced to study characteristic classes of
bundles over orbifold loop spaces. For our statement to be true as stated, we
will need to invert $q$. {\em Twisted} here means
that we will have to replace the centralizers with the appropriate central
extensions and $q^\frac{1}{|g|}$ with $q^\frac{1}{h|g|}$, where $h$ is
the order of 
$\alpha\at{B\gg}$.
The equivariant Witten genus takes 
values here. 
In \cite{Ganter:stringypowers}, I defined power operations $P_n$ in 
an appropriate subring of Devoto
equivariant Tate $K$-theory. If a function $f(g,h;\tau)$ satisfies
condition (a) of Norton's conjecture, this implies that $f$ is in this
subring. 
My definitions in \cite{Ganter:stringypowers} were guided by
the work of Dijkgraaf, Moore, Verlinde and Verlinde on orbifold genera
of symmetric powers. The results of \cite{Ganter:stringypowers},
together with \cite{Devoto:elliptic}, make
$K_{\Dev,G}$ an equivariant theory with power operations and a
Hopkins-Kuhn-Ravenel character theory. In this situation, 
we have the following combinatorial definition of Hecke operators
(cf.\ \cite{Ando:thesis}, and \cite{Ganter:thesis}):
$$
  T_n(f) =
  \frac1n\sum_{\stackrel{[\sigma,\rho]}{\text{transitive}}} P_n(f)(\sigma,\rho), 
$$
where the sum runs over all conjugacy classes of pairs of commuting elements
$(\sigma,\rho)$ of the symmetric group $\Sn$ with the property that 
$\langle\sigma,\rho\rangle$ acts transitively on $\{1,\dots,n\}$.

The main result of this paper can be summarized as
follows: 
\begin{Thm}
  The following three definitions of Hecke operators result in exactly
  the same formula for functions $f(g,h;\tau)$:
  \begin{enumerate}
  \item The geometric definition using isogenies, when $f$ is an
    element of (twisted) Ginzburg-Kapranov-Vasserot equivariant elliptic cohomology;
  \item The combinatorial definition, using the
    equivariant power operations of \cite{Ganter:stringypowers}, when
    $f$ is an element of Devoto's equivariant Tate $K$-theory on which
    $P_n$ is defined;
  \item for $g=1$, the definition of twisted Hecke operators acting on
    McKay-Thompson series.
  \end{enumerate}
  The formula is 
  $$
    T_n(f)(g,h;\tau) = \sum_{\stackrel{ad=n}{0\leq b<d}}
    f\(g^d,g^{-b}h^a;\frac{a\tau +b}{d}\).  
  $$
\end{Thm}
Similarities between equivariant elliptic cohomology and generalized
Moonshine have been studied in the past by Baker, Devoto, and Thomas
(cf.\ e.g.\ \cite{Baker:Thomas} and \cite{Devoto:elliptic}). Baker's
work includes a definition of twisted Hecke 
operators in elliptic cohomology, using Adams operations. The theory of
Hecke operators in elliptic cohomology was developed further in Ando's
work, which brought the theory of isogenies of formal groups
into the picture \cite{Ando:poweroperations}. 
\subsection{Acknowledgements}
It is a pleasure to dedicate this paper to John McKay, who introduced
me to the subject of 
generalized Moonshine. 
Many thanks go to Jorge Devoto, from whom I learned a lot
about the connections between Moonshine and elliptic cohomology, to
Jack Morava and Dan Freed for their very helpful email
correspondence, to Alex Ghitza for plenty of help and comments on earlier
versions of the paper, to Michael Hopkins, Matthew Ando, and Christian
Haesemeyer for long and helful conversations about the subject.
\section{Principal bundles over complex elliptic curves}\label{Prin-ell-Sec}
Let $G$ be a finite group, let $\mC(G)$ denote the set of all pairs of
commuting elements in $G$, and consider
the right action of $G$ on $\mC(G)$ given by simultaneous conjugation
$$
  (g,h)\mapsto (s\inv gs,s\inv hs).
$$ 
We define
the {\em generalized centralizer} group
$C_G(g,h)$ of the pair $(g,h)$
by $$C_G(g,h)=\Stab_G(g,h).$$
The {\em generalized conjugacy class} $[g,h]_G$ is defined as the orbit of
$(g,h)$ under $G$.
The set of all conjugacy classes of pairs of commuting elements in
$G$ is denoted $\PG$. 

Let $E$ be an elliptic curve
over $\CC$. 
Assume that we have picked an isomorphism
\begin{equation}
  \label{Lt-Eqn}
  E\cong \CC/\langle\tau,1\rangle.
\end{equation}
Throughout the paper, we will identify $\Lt$ with $\ZxZ$ via the
isomorphism
\begin{eqnarray}\label{ZxZ-to-Lt-Eqn}
\notag  \ZxZ&\to&\Lt\\
  (1,0)&\mapsto& -1\\
\notag  (0,1)&\mapsto& \tau
\end{eqnarray}
Further, we will view $\CC$ with the action of the additive subgroup
$\langle\tau,1\rangle\cong\ZZ^2 $ 
as a contractible free $\ZZ^2$-space
$E\ZZ^2$, and $E$ as a classifying space 
$$
  E=\BZZ=E\ZZ^2/\ZZ^2. 
$$
\begin{Lem}
  The set of isomorphism classes of principal $G$-bundles over $E$ is
  $$
    \Prin GE\cong [\BZZ,\BG]\cong \Hom(\ZZ^2,G)/\on{Inn}(G) =\PG,
  $$
  the set of all conjugacy classes of pairs of commuting elements of $G$.
\end{Lem}
Such a pair of commuting elements
may also be viewed as encoding the monodromy of the corresponding
principal bundle along the
two standard circles of the torus.

Under the isomorphism (\ref{ZxZ-to-Lt-Eqn}), the map
\begin{eqnarray*}
  \varphi\negmedspace :\Ltp & \longrightarrow & \Lt \\
\notag  \tau' & \mapsto & a\tau + b \\
\notag   1     & \mapsto & c\tau + d
\end{eqnarray*}
corresponds to the endomorphism of $\ZZ^2$ given by the matrix
\begin{equation}\label{phi-Eqn}
  \varphi = \(
  \begin{array}{rr}
    d&-b\\
    -c&a
  \end{array}
  \),
\end{equation}
and
\begin{eqnarray*}
  \varphi\dual\negmedspace :\Lt & \longrightarrow & \Ltp \\
  \tau & \mapsto & d\tau' -b \\
  1     & \mapsto & -c\tau' + a
\end{eqnarray*}
becomes the pseudo-inverse of this matrix,
\begin{equation*}
  \varphi\dual = \(
  \begin{array}{rr}
    a&b\\
    c&d
  \end{array}
  \).
\end{equation*}

The map $\phi :=B\varphi$ is an isogeny
$$\map{\phi}{\CC/\Ltp}{\CC/\Lt}$$ with
dual isogeny $\phi\dual=B(\varphi\dual)$.
Any isogeny of complex elliptic curves can be obtained in this way. If
$ad-bc=1$, then $\phi$ becomes an isomorphism with $\phi\inv=
\phi\dual$. It is the isomorphism obtained when comparing two
different choices for (\ref{Lt-Eqn}).

\begin{Prop}
  The isogeny $\phi\dual$ pulls back the principal bundle over
  $\CC/\Ltp$ which is classified by the pair $(g,h)$ to the bundle over
  $\CC/\Lt$ which is 
  classified by 
  $(g^ah^c,g^bh^d)$.
\end{Prop}
\begin{Pf}{}
  After identifying both $\Lt$ and $\Ltp$ with $\ZxZ$ as in
  (\ref{ZxZ-to-Lt-Eqn}),  
  the map
  $\varphi\dual$ sends $(1,0)$ to $(a,c)$ and 
  $(0,1)$ to $(b,d)$. Recall that $(g,h)$ stands for the map from
  $\ZxZ$ to $G$ sending $(1,0)$ to $g$ and $(0,1)$ to $h$. If we
  precompose this map with $\varphi\dual$, we obtain
  $$(g,h)\circ\varphi\dual = (g^ah^c,g^bh^d).$$
\end{Pf}
\subsection{The moduli space of principal $G$-bundles}\label{moduli-space-Sec}
Let
$$
(s,\gamma)\in G \times \SL
$$
act on
$$
  \mC(G)\times\HH
$$
from the right by 
$$
  (g,h;\tau)\mapsto \(s\inv(g^ah^c,g^bh^d)s;\gamma^{-1}(\tau)\),
$$
where
$$
  \gamma =  \left(
    \begin{array}{rr}
      a&b \\
      c&d 
    \end{array}\right) \quad\text{and}\quad\gamma\inv(\tau) =
  \frac{d\tau-b}{-c\tau+a} . 
$$
Then the quotient 
$$
  \PG\times\HH \medspace \slash \medspace\SL
$$
is the coarse moduli space $X_G$ 
of isomorphism classes of principal $G$-bundles over complex elliptic curves.
Apart from the roots of unity, conditions (a)
and (c) of Norton's Conjecture are precisely the coherence data
satisfied by
holomorphic functions on $X_G$ (compare also \cite{Devoto:elliptic}). 

In order to understand the roots of unity, we need to consider the moduli stack
$$
  \mM_G := \mC(G)\times\HH \medspace\slash\!\!\slash\medspace G\times \SL .
$$
Write
$$
  S:= \left(
    \begin{array}{rr}
      0&1 \\
      -1&0 
    \end{array}\right)\quad\text{ and }\quad 
  T:=\left(
    \begin{array}{rr}
      1&1 \\
      0&1 
    \end{array}\right).
$$
%
The groupoid $\mC(G)\mmod G$ is equivalent to the category of 
principal $G$-bundles over the torus and isomorphisms between them. 
To show this, consider the pointed circle $(\bbS^1,*)$ and the torus
$$(\bbS^1,*)\times(\bbS^1,*).$$  
Let $(P,p)$ be a principal $G$-bundle over this torus together with a
choice of basepoint over $(*,*)$. 
Let $(g,h)$ be the commuting pair
given by the monodromy of $P$ around the two circles, starting at $p$.
Note that for any two pointed principal bundles $(P,p)$ and $(P',p')$
with the same monodromy, there is a unique basepoint-preserving 
isomorphism from $P$ to $P'$.

For any commuting pair $(g,h)$ in $G$, fix once and for all a pointed
principal $G$-bundle $(P_{g,h},p)$ over the torus with this monodromy.
A (not necessarily basepoint-preserving) morphism from $(P,p)$ to
$(P',p')$ is a pair $(f,s)$, 
where $s\in G$, and $\map f{(P,ps)}P'$ is a basepoint-preserving
isomorphism. For any $s\in G$, the pointed bundle $(P_{g,h},ps)$ has
monodromy $(s\inv gs,s\inv hs)$, and 
there is a unique basepoint-preserving
isomorphism  
$$
  \map f{(P_{g,h},ps)}{(P_{s\inv gs,s\inv hs},p')}.
$$
%

Note that the automorphism group of $(g,h)$ in $\mC(G)\mmod G$ is
given by the
generalized centralizer group $C_G(g,h)$. The automorphisms of
$(g,h;\tau)$ in $\mM_G$ might be more complicated: 
The automorphisms of $P_{g,h}$ over $E$ are parametrized by
$C_G(g,h)$, but now 
there might also be automorphisms of $E$ which pull back $P_{g,h}$ to a
bundle isomorphic to $P_{g,h}$. 
For instance, $(1,ST)$
is an automorphism of $\(g,g\inv;\frac{\sqrt3i-1}{2}\)$ whenever $g$ has
order three, and $(s,S)$ is an automorphism of $(g,h;i)$ whenever
$s(g,h)s\inv = (h\inv,g)$. 
\subsection{Line bundles over $\mM_G$}\label{Freed-Quinn-Sec}
\subsubsection{The cocycle $\alpha$}
Fix, once and for all, a cellular decomposition of $BG$. 
Let $\map{i_{(2)}}{BG^{(2)}}{BG}$ denote the inclusion of the $2$-skeleton
in $BG$. 
\begin{Def}[{compare \cite{Willerton}}]
  A cocycle $\alpha\in Z^3(BG;\mathbb R/\ZZ)$ is called {\em
    normalized}, if $i_{(2)}^*(\alpha)=0$. 
\end{Def}
\begin{Lem}
  Any cocycle $\alpha\in Z^3(BG;\mathbb R/\ZZ)$ is cohomologous to
  a normalized one.
\end{Lem}
\begin{Pf}{}
  The statement follows from the commuting diagram
  $$
    \xymatrix{
     {C^2(BG;\mathbb R/\ZZ)}\ar@{->>}[r]^{i_{(2)}^*{\phantom {xx}}}\ar[d]_\delta & 
                             {C^2(BG^{(2)};\mathbb R/\ZZ)}\ar@{->>}[d]^{\delta} \\
     {Z^3(BG;\mathbb R/\ZZ)}\ar[r]^{i_{(2)}^*{\phantom {xx}}} & 
                             {Z^3(BG^{(2)};\mathbb R/\ZZ).}
    }
  $$
  Here $C^*$ are the groups of cochains and $Z^*$ the groups of
  cocycles. 
  The upper horizontal arrow is surjective, since $\mathbb R/\ZZ$ is
  an injective abelian group.  
\end{Pf}
We fix a normalized cocycle $\alpha\in
Z^3(BG;\mathbb R/\ZZ)$. 
Recall that we have
$$H^4(BG;\ZZ)\cong H^3(BG;\mathbb R/\mathbb Z).$$
In the case that $G$ is the Monster, a candidate for $[\alpha]\in
H^4(BG;\ZZ)$ was described in \cite{Mason:cohomology}. 
\subsection{Construction and properties of the Freed-Quinn line bundle $\mL^\alpha$}
Inspired by the work of Dijkgraaf, Vafa, Verlinde and Verlinde \cite{DVVV}, 
Freed and Quinn use $\alpha$ to define a (possibly degenerate) line bundle
$L^\alpha = \{L_{g,h}^\alpha\}$ on $\overline\mC(G)$. 
%
They proceed to describe the $\SL$-action on this line bundle. 
Their discussion can be interpreted as constructing an $\SL$-equivariant
line bundle
$$
  \xymatrix{
  {\mathfrak{H}\times L^\alpha} 
  \ar[d]\\
  \HH\times\overline\mC(G),}
$$
and thus a line bundle $\mL^\alpha$ on the moduli stack
$\mM_G$.  
We briefly recall their construction, which simplifies when $\alpha$
is normalized.
\begin{Cons}\label{L-Con}
  Let $P$ be a principal $G$ bundle over $\mathbb T^2$, and consider
  the groupoid $\mC_P$ whose objects are cellular classifying
  maps $\map fPEG$ for $P$, and whose morphisms from $f$ to $f'$
  are homotopy classes
  rel boundary of cellular $G$-homotopies
  $$
    \map f{[0,1]\times P}{EG}
  $$
  from $f$ to $f'$. We define a functor $\mF$ from $\mC_P$ to the
  category of metrized complex lines as follows:
  For every $f$, one sets\footnote{The metrized integration line
  $I_{\mathbb T^2,\overline{f}^*\alpha}$ of \cite{Freed:Quinn} is
  canonically trivialized if 
  $\alpha$ is normalized and $f$ is cellular.} 
  $\mF(f)=\CC$.
  Further $\mF(f\stackrel
  h\to f')$ is multiplication with $$e^{2\pi i(\overline
  h^*\alpha)(x)},$$ where $x$ is a $3$-cycle representing the
  fundamental class in $[0,1]\times \TT^2$, and $\overline h=h/G$.
  It can be shown that this is independent of the choice of
  representatives $h$ and $x$ and that $\mF$ has no holonomy. 
  The line $L^\alpha_P$ is then defined as the space of invariant
  sections of $\mF$.

  Let $\psi$ be an automorphism of $P$ (not necessarily covering
  the identity on $\mathbb T^2$). Then the action of $\psi$ on
  $L^\alpha_P$ is given by multipication with the {\em Chern-Simons}
  invariant of the glued mapping cylinder 
  $$
    \cyl\psi P:=[0,1]\times P\medspace\slash \sim\medspace, \quad
    \text{where}\quad(0,\psi(p)) \sim (1,p).
  $$
  This invariant is computed by choosing a cellular classifying map
  $\map h{\cyl\psi P}EG$ and a $3$-cycle $x$ representing the fundamental
  class of $\cyl{\overline{\psi}}{\TT^2}$ and then evaluating $\overline
  h^*(\alpha)$ at $x$. More generally,
  let $\gamma\in\SL$ be an automorphism
  of $\mathbb T^2$. Then the induced isomorphism $\gamma_*$ of metrized
  integration lines is defined to be the identity 
  \begin{equation}
    \label{gamma*-Eqn}
    \map{\id_\CC}{\mathcal F(f)=\CC}{\CC=\mathcal
    F(f\circ\gamma)}.
  \end{equation}
\end{Cons}
Note that $\TT^2$ can be replaced with any Riemann surface.
We will often drop $\alpha$ from the notation and write $L_P$ for
$L_P^\alpha$. If $P$ is a principal bundle over $\TT^2$ classified by
$(g,h)$, we will sometimes write $L_{g,h}$ for $L_P$.
\begin{Prop}\label{Ln-Prop}
  Let $\map\phi{E'}E$ be an isogeny of degree
  $n$, and let $P$ be a principal $G$-bundle over $E$. Then there is
  an isometry of metrized lines, 
  $$
    \map{\phi_*}{L_{\phi^*(P)}}{L_P^{\tensor n}}.
  $$
  The construction of $\phi_*$ is natural in $E$.
\end{Prop}
\begin{Pf}{} Write $\mathbb I$ for $[0,1]$.
  Let $A$ be the kernel of $\phi$, and let $C^A_3(\mathbb I\times E')$ denote the
  group of $A$-invariant singular $3$-chains in $\mathbb I\times E'$. Since $A$ acts
  fixed point free on $\mathbb I\times E'$, every element $x\in C^A_3(\mathbb I\times
  E')$ is of the form  
  $$
    x = \sum \lambda_i \(\sum_{a\in A} a\)s_i,
  $$
  where the $\lambda_i$ are integers and the $s_i$ are singular
  simplices in $\mathbb I\times E'$. Hence the map $\id_\mathbb I\times\phi_*$ 
  takes values in $n C_3(\mathbb I\times E)$.
  If $x\in C_3^A(\mathbb I\times E')$ represents the relative fundamental
  class of $\mathbb I\times E'$, then  
  $\frac1n\phi_*(x)$ represents the relative fundamental class of
  $\mathbb I\times E$. 
  Let $\mC_P$ and  $\mC_{\phi^*P}$ be the groupoids in Construction
  \ref{L-Con}. Precomposition with $\phi$ gives a map of groupoids
  $$
    \map\iota{\mC_P}{\mC_{\phi^*P}}.
  $$
  Let $\mF_{P,\alpha}$ and $\mF_{\phi^*P,\alpha}$ be the functors in
  Construction \ref{L-Con}. 
  The discussion above gives rise to a natural isomorphism
  $$
    \mF_{P',\alpha}\circ\iota\cong\(\mF_{P,\alpha}\)^{\tensor n}.
  $$
  It remains to observe that the invariant sections of
  $\mF_{P',\alpha}$ and $\mF_{P',\alpha}\circ\iota$ are canonically
  isomorphic. 
\end{Pf}
For later reference, we note that the definition of the isomorphism in
Proposition \ref{Ln-Prop} used only the structure of $\phi$ as a
principal $A$-bundle. Note also that the action of $A$ on
$L_{\phi^*P}$ is trivial. 
\begin{Lem}\label{alpha=0-Lem}
  Let $\alpha=0$. Then the
  line bundle $L^\alpha$ over $\overline{\mC}(G)$ is canonically
  trivial, and the $\SL$-action preserves this trivialization. 
\end{Lem}
\begin{Pf}{}
  Let $\alpha=0$, then the line of the previous construction is
  canonically trivialized. Explicitly, such a trivialization is given
  by picking an arbitrary object $f\in \ob(\mC_P)$ and evaluating invariant
  sections of $\mF$ at $f$. 
  Since our choice of trivialization 
  was independent of the choice of $f$, it follows that
  $$\map{\gamma_*=\id_\CC}{L^\alpha_P}{L^\alpha_{\gamma^*P}}.$$
\end{Pf}
\begin{Lem}\label{change-of-groups-Lem}
  Let $\map aHG$ be a map of finite groups. Then $a$ induces a map of
  groupoids
  $$\map{\mC(a)}{\mC(H)\mmod H}{\mC(G)\mmod G},$$ which
  pulls back
  $L^{\alpha}$ to
  $L^{a^*(\alpha)}$ as an $\SL$-equivariant bundle. 
\end{Lem}
\begin{Pf}{}
  Fix a cellular classifying map 
  $$
    \map{F}{EH[G]}{EG},
  $$
  and write $Ba$ for the induced map $\overline F$ of classifying
  spaces. Let $P$ be a principal $H$-bundle over $\mathbb T^2$, 
  let $\mathcal C_P$ be the groupoid in \cite[p.6]{Freed:Quinn}, and
  similarly, $C_{P[G]}$. We have an equivalence of
  groupoids 
  \begin{eqnarray*}
       \iota\negmedspace:
    \mC_P&\to&\mC_{P[G]}\\
       f &\mapsto& F\circ f[G]\\
       h &\mapsto& F\circ h[G].
  \end{eqnarray*}
  Here $\map fPEH$ ranges over the classifying maps of $P$, and $h$ ranges over
  the homotopies between such classifying maps.
  Let $\mathcal F_{P,Ba^*(\alpha)}$ be the functor defined in
  \cite[p.6]{Freed:Quinn}, and note that 
  $$\mathcal F_{P,Ba^*(\alpha)}=\mathcal
  F_{P[G],\alpha}\circ\iota.$$
  It follows that the corresponding metrized lines of invariant
  sections of the two functors are the same:
  $$
    L_{P}^{Ba{\phantom{\vert\negmedspace}}^*(\alpha)} = L_{P[G]}^{\alpha}.
  $$
  Let now $F'$ denote a different choice for $F$ and let $\iota'$ be
  the corresponding equivalence of groupoids. Then there exists a
  cellular 
  homotopy between $F$ and $F'$ and hence a natural transformation
  between $\iota$ and $\iota'$, yielding a canonical isomorphism
  between the spaces of invariant sections of $\iota^*\mathcal F$ and
  $(\iota')^*\mathcal F$. The construction of the proof is compatible
  with the functoriality condition in the following sense: Let
  $\map\psi PP'$ be a map of principal $H$-bundles covering an
  orientation-preserving diffeomorphism of the base (for the proof of
  this lemma it
  is enough if this isomorphism is the identity of the torus). Then
  the induced isometries $\psi_*$ and $\psi[G]_*$ of the above lines agree.
\end{Pf}
\begin{Cor}\label{L1-Cor}
  Let $G$ be a finite group, and let $\alpha$ be
  a normalized $3$-cocycle on $BG$ with values in $\mathbb R/\ZZ$. 
  Then the line
  $L_{1,1}^{\alpha}$ is canonically trivialized and carries a
  trivial $\SL$-action.
\end{Cor}
\begin{Pf}{}
  By the previous lemma, it is enough to consider the case where $G=1$
  is the trivial group. Let $\map{\id_{\mathbb T^2}}P{\mathbb T^2}$ denote
    the trivial principal bundle over $\mathbb T^2$.
  Since trivial bundles can be classified by constant maps,
  we have an equivalence of groupoids
  $\map\iota{\mC_{P}^{(2)}}{\mC_{P}}$, where
  $\mC_{P}^{(2)}$ is the subgroupoid whose objects are
  classifying maps $\map f{\mathbb T^2}{B1^{(2)}}$ and whose morphisms are
  homotopy classes relative boundary of homotopies 
  $\map h{[0,1]\times\mathbb T^2}{B1^{(2)}}$.
  Since $\alpha$ is normalized, it follows that
  $\iota^*(L^\alpha_P)$, and hence $L^\alpha_P$ is canonically
  trivialized. 
\end{Pf}
\begin{Cor}
  Let $P_{g,h}$ be a principal bundle over the torus which is
  classified by $[g,h]_G$. Let $n=|g|\cdot|h|$. Then 
  $\(L_{[g,h]}^\alpha\)^{\tensor n}$ is canonically trivial.
\end{Cor}
\begin{Pf}{}
  The trivialization is given by Proposition \ref{Ln-Prop} applied to
  the $n$-fold covering which wraps the first circle around itself
  $|g|$ times and the second one around itself $|h|$ times. 
  This covering pulls back $P_{[g,h]}$ to a trivial bundle.
\end{Pf}

\begin{Lem}
  Let $p$ be a map from the torus to the circle, and let
  $P$ be the pull-back along $p$ of a principal $G$-bundle $Q$
  over the pointed circle. 
  Then $p$ induces a trivialization of $L^\alpha_{P}$.
\end{Lem}
\begin{Pf}{}
  Let $f$ be a classifying map for $Q$. Then $f\circ p$ is a
  classifying map for $P$ and hence yields a trivialization of
  $L_P^\alpha$. If $f'$ is a different classifying map for $Q$, then there
  is a cellular homotopy $\map h{\mathbb I\times\bbS^1}{EG}$, and
   $\overline{h}\circ(\id_{\mathbb I}\times p)^*(\alpha)=0$. Hence our
   trivialization is 
  independent of our choice of $f$.
\end{Pf}
For every $g\in G$, we fix, once and for all, a principal $G$-bundle
$P_g$ 
over the pointed circle with monodromy $g$.
We will follow the convention to trivialize lines of the form
$L_{g,1}^\alpha$ or $L_{1,g}^\alpha$ using the projection to the first
(respectively second) generating circle of the torus.
Let $d$ be a natural number, and let $f_g$ be a classifying map for $P_g$. 
Precomposing $f_g$ 
with the degree $d$ self-map of the circle gives a classifying map for
$P_{g^d}$. This yields isomorphisms of
trivialized lines 
$$
  i_{d}\negmedspace : L_{1,g^d}\cong \CC=\CC \cong L_{1,g},
$$
which are compatible in the sense that $i_{ad}=i_d\circ i_a$ and
$i_{k|g|+1}=\id_{L_{1,g}}$. 
We will occasionally use this isomorphism to identify these lines and
write $L_\gg:=L_{1,g}$.
Here $\gg$ stands for the cyclic subgroup generated by $g$.
Let 
$$
\gamma =    \(\begin{array}[]{cc}
    a&b\\ 
    nc&d
    \end{array}\)\in\Gamma_0(n).
$$
Then $ad\equiv 1\mod n$, and $\langle g^d\rangle=\gg$. Hence 
$$
  \map{\gamma_*^{-1}}{L_{1,g}}{L_{1,g^d}}
$$
is an automorphism of $L_\gg$. 
We caution the reader that this action of $\Gamma_0(n)$ on $L_\gg$ 
is not independent of the
choice of generator $g$. Instead, replacing $g$ by $g^l$ amounts to
pulling back $\alpha\at{B\gg}$ along the map $(-)^l$.
We will see in the proof of Lemma \ref{restriction-Lem} below that
restriction along this map raises the action to the $l^{th}$ power.
\begin{Lem}\label{Gamma-0-Lem}
Let $h=|\alpha\at{B\gg}|$. Then 
$T^n$
acts as multiplication by an $h^{th}$ root of unity 
on the line $L^\alpha_{g,1}$. 
Similarly, any element
$\gamma\in\Gamma_0(n)$
acts as multiplication by a $h^{th}$ root of unity
on the line $L^\alpha_\gg$.
\end{Lem}
We will see in Corollary \ref{h-Cor} below, that $T^n$ actually acts by
multiplication with a primitive $h^{th}$ root of $1$.
\begin{Pf}{}
  Let $\phi$ be the $n$-fold covering map from the torus to itself
  which
  has degree $n$ on the first circle and is the identity on the
  second.
  Then $\phi$ pulls back $P_{g,1}$ to the trivial bundle, and we have
  $$
    T^n\circ\phi=\phi\circ T.
  $$
  Hence $\phi$
  induces an $n$-fold covering of the glued cylinders of Construction
  \ref{L-Con}:  
  $$
    \map\phi{\cyl T{P_{1,1}}}{\cyl{T^n}{P_{g,1}}}.
  $$
  Just as in the proof of Proposition \ref{Ln-Prop},
  we may conclude that the Chern-Simons invariant of ${\cyl
    T{P_{1,1}}}$ is equal to the $n^{th}$ power of the Chern-Simons
  invariant of 
  ${\cyl{T^n}{P_{g,1}}}$. The first gives the action of $T$ on
  $L_{1,1}^\alpha$. By Corollary \ref{L1-Cor}, this action is $1$.
  The Chern-Simons invariant of   ${\cyl{T^n}{P_{g,1}}}$ gives the action of $T^n$ on
  $L_{g,1}^\alpha$. It follows that this action is an $n^{th}$ root of
  $1$. From this point of view, it is also clear why its order
  divides the order of $\alpha$.
Fix $g$, and let $\gamma\in\Gamma_0(n)$ be as above. The same
argument, using the degree $n$ map of the second circle, shows that
for any two classifying maps $f$ and $f'$ of $P_{1,g^d}$, the
corresponding trivializations of $L_{1,g^d}$ differ by an $n^{th}$ root
of unity. Writing $f_g$ and $f_{g^d}$ for the classifying maps used to
trivialize $L_{1,g}$ and $L_{1,g^d}$, the claim follows when we take
$f=f_{g^d}$ and $f'=f_{g}\circ\gamma\inv$. 
\end{Pf}
\subsection{Sections of $L^\alpha$}\label{Sections-Sec}
The following definitions and facts are taken from
\cite[5]{Freed:Quinn}. 
We let $E$ be the space of sections of $L^\alpha$ over
$\overline\mC(G)$. 
In the notation of \cite{Freed:Quinn}, 
we have $$E=\mathcal{L}^2(\overline{\mathcal{C}}_{\mathbb{S}^1\times\mathbb{S}^1},
  \overline{\mathcal{L}}_{\mathbb{S}^1\times\mathbb{S}^1}).$$
Further, Freed and Quinn define a line bundle
$\overline{\mathcal{L}}_{[0,1]\times\mathbb{S}^1}$ on the set
$\overline{\mC'}_{[0,1]\times\mathbb{S}^1}$  
of isomorphism classes of principal $G$-bundles over
$[0,1]\times\mathbb{S}^1$ with  basepoints chosen over the basepoints
of the boundary. They prove that the sections of this line bundle
$$
  A:=\mathcal{L}^2(\overline{\mathcal{C}'}_{[0,1]\times\mathbb{S}^1},
  \overline{\mathcal{L}}_{[0,1]\times\mathbb{S}^1})
$$
form a coalgebra, whose underlying vector space is the direct sum of
complex lines
$$
  A \cong \bigoplus_{g,h\in G}L_{g,h}.
$$
If we view $A$ as the space of complex valued linear functions on its dual
algebra 
$$
  A^* \cong \bigoplus_{g,h\in G}L^*_{g,h},
$$
then $E\sub A$ is identified with the subspace of central
  sections.

The algebra $A^*$ is semisimple, and $E$ possesses an orthonormal
basis of character functions $\{\chi_\lambda\}$, where $\lambda$ runs
through the irreducible representations of $A^*$. The support of
$\chi_\lambda$ is 
$$
  \on{supp}(\chi_\lambda) = \{(g,h)\mid g\in[g_0],h\in C_g\}
$$
for a conjugacy class $[g_0]$ depending on $\lambda$. For fixed $g$,
the restriction of $\chi_\lambda$ to (the units of)  
$$
  \bigoplus_{h\in C_{g}}L^*_{g,h}
$$
is a character of the central extension $\hat C_g$ of $C_g$ defined by
these units.
Thus 
$$
  E \cong \CC\{\chi_\lambda\}\cong\bigoplus_{[g]\sub G}R(\hat C_g)\tensor\CC.
$$
Following Freed and Quinn, we write 
$e_{g,1}\in L_{g,1}$ and $e^*_{g,1}$ for the elements of unit length
corresponding to the trivialization we picked above and
note that $e^*_{g,1}$ is the unit of the group $\widehat C_g$.
\subsubsection{Sections of $\mL^\alpha$}
By construction, a section of $\mL^\alpha$ is equivalent to an
$\SL$-equivariant map
$$
  \map{f}{\HH}{E.}
$$
Fix an element $g\in G$ and consider the $[g]$-component of $f$,
$$
  \map{f_g}{\HH}{R(\hat C_g)\tensor\CC.}
$$
\begin{Lem} The element
$T^n\in\SL$
acts as multiplication by a root of unity 
on the summand $R(\hat C_g)\tensor\CC$. Its order divides $n$ and the
order of $\alpha$. More precisely $T^n$ acts by
multiplication with $e^{2\pi i\alpha(x)}$, 
where $x$ is the $3$-cycle 
$$x=\sum_{k=0}^n(g,g^k,g)\in Z_3(B\gg).$$
\end{Lem}
\begin{Pf}{}
  \begin{figure}[h]
    \centering
\includegraphics[angle=88,scale=.5292]{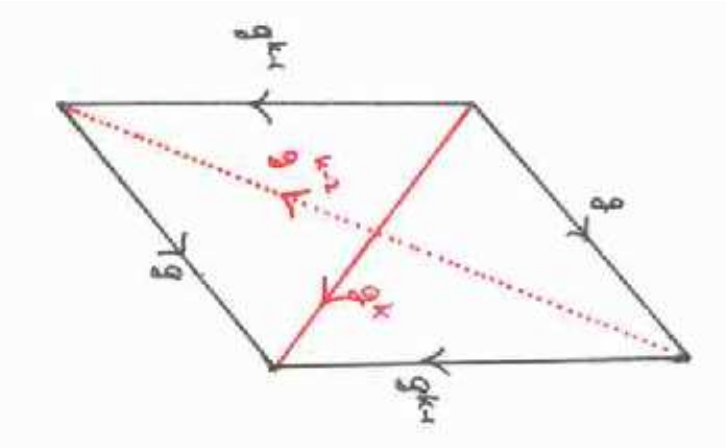}
\hspace{-.5cm}
\includegraphics[angle=90,scale=0.7]{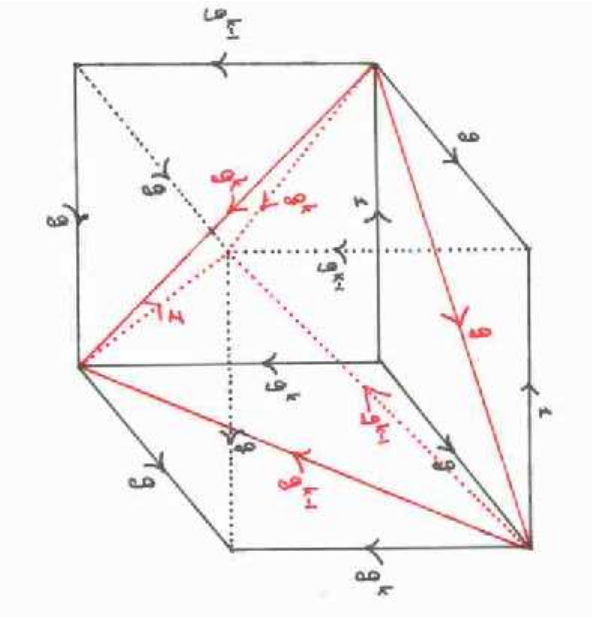}
\hspace{-1cm}
\includegraphics[angle=90,scale=0.5292]{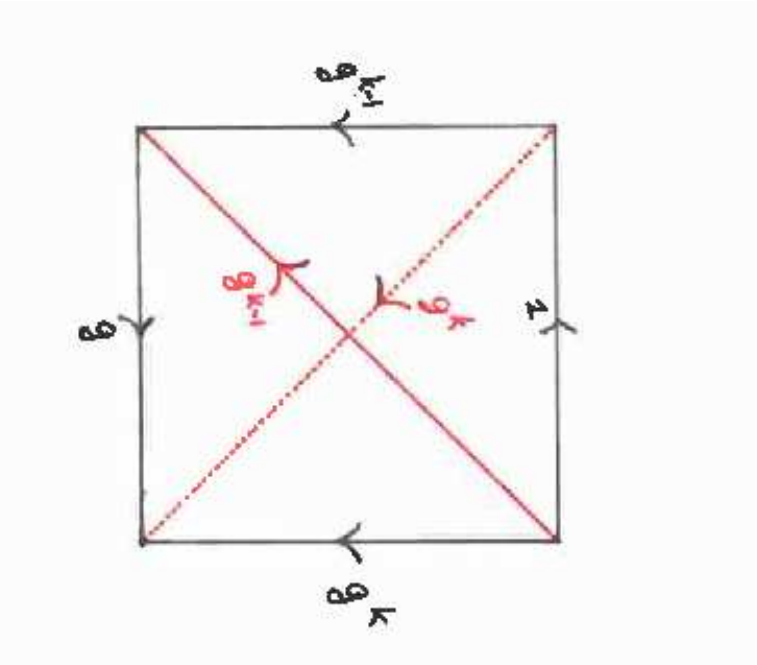}
    \caption{The triangulation of the cube (middle) and the degenerate
      simplices.}
\label{cube:Fig}
  \end{figure}
  The action of $T$ on $A^*\at{[g]}$ is given by multiplication
  with the element $e^*_{g,g}\in L^*_{g,g}$ defined in \cite[(5.5)]{Freed:Quinn}.  
  Hence it is enough to understand the action of $T^n$ on
  $L_{g,1}^\alpha$. We saw in Lemma \ref{Gamma-0-Lem} that this action
  is given by multiplication with a root of unity whose order divides
  $n$ and the order of $\alpha$. To determine it explicitly, we
  triangulate the fundamental domain $\mathbb I\times \mathbb I$ of
  $\TT^2=\mathbb R/\ZZ$ by choosing as $1$-simplices $s_1$ from $(0,0)$
  to $(0,1)$, $s_2$ from $(0,0)$ to $(1,0)$, $s_3$ from $(0,1)$ to
  $(1,1)$, $s_4$ from $(1,0)$ to $(1,1)$ and $s_5$ from $(0,1)$ to
  $(1,0)$. Then $T$ pulls back the corresponding triangulation of
  $\TT^2$ to the one corresponding to the following choice of
  $1$-simplices on $\mathbb I^2$: a simplex $t_1$ from $(0,1)$ to
  $(0,0)$, $t_2$ from $(0,0)$ to $(1,0)$, $t_3$ from $(0,1)$ to
  $(1,1)$, $t_4$ from $(1,1)$ to $(1,0)$, and $t_5$ from $(0,0)$ to
  $(1,1)$. Cutting open the pointed bundle $P_{g,g}$ along the two
  standard circles of $\mathbb T^2$ gives a bundle over $\mathbb I^2$
  with basepoints over the corners and parallel transports $g$ along
  $s_1$, $s_2$, $s_3$ and $s_4$, and $1$ along $s_5$. The same
  procedure for $T^*P_{g,g}=P_{g,1}$ gives parallel transports 
  $g$ along $t_2$, $t_3$, and $t_5$, and $1$ along $t_1$ and $t_4$. We
  triangulate the faces of the cube $\mathbb I^3$ by choosing the
  first of the above triangulations on $\mathbb I^2\times\{0\}$ and
  the second one on $\mathbb I^2\times\{1\}$ and adding the
  $1$-simplices $r_1$ from $(1,1,1)$ to $(0,1,0)$, $r_2$ from $(1,1,1)$ to
  $(1,0,0)$, $r_3$ from $(0,0,1)$ to $(0,1,0)$, $r_4$ from $(0,0,1)$
  to $(1,0,0)$, $r_5$ from $(1,1,1)$ to $(1,1,0)$, $r_6$ from
  $(1,0,1)$ to $(1,0,0)$, $r_7$ from $(0,1,1)$ to $(0,1,0)$, and $r_8$
  from $(0,0,1)$ to $(0,0,0)$. Note that there exists a bundle $Q_k$
  over $\mathbb I^3$, pointed over the corners, such that the parallel
  transports along $r_5$ and $r_7$ are $g^{k-1}$ and those along $r_6$ and
  $r_8$ are $g^{k}$. The cube $\mathbb I^3$ is divided into five
  $3$-simplices whose edges are the $1$-simplices specified above, and
  the cycles defined by pushing forward these $3$-simplices along a
  classifying map of $Q_k$ are
  $$
    \pm(g^{k-1},g,1),\quad \pm(g^{k-1},1,g),\quad \pm(1,g,g^{k-1}),\quad
    \pm(g,1,g^{k-1})\quad\text{and} \quad\pm(g,g^{k-1},1).  
  $$ 
  We add two more degenerate $3$-simplices, one to the face
  $\{0\}\times\mathbb I^2$ giving the cocycle $\pm(1,g^{k-1},g)$ and one
  along the face $\mathbb I\times\{0\}\times\mathbb I$ giving the
  cocycle $\pm(g,g^{k-2},g)$.
  This allows us to glue the first circle of the torus back
  together via 
  $$
    (0,x,y)\sim(1,x,y).
  $$ 
  This gives a triangulation of
  the mapping cylinder of $T\at{\bbS^1\times\mathbb I}$.
  Further, we are able to glue $Q_k$ to $Q_{k+1}$ along the ends of
  the cylinder in such a way that their triangulations fit together.
  Since $Q_{n} = Q_0$, we may join $Q_0$ through $Q_{n-1}$ together to
  get a bundle over the mapping cylinder of $T^n$.
  Since $\alpha$ is normalized, it vanishes on cycles containing the
  identity. 
  Hence $T^n$ acts as claimed. 
\end{Pf}
\begin{Lem}\label{generator-Lem}
  The cycle $\sum_{k=0}^{n-1}(g,g^k,g)$ is a generator of the cyclic
  group $$H_3(B\gg)\cong\ZZ/n\ZZ.$$
\end{Lem}
\begin{Pf}{}
  We will use the terminology and notation of \cite{Brown}. 
  Let $\ZZ\leftarrow F_\bullet$ be the standard resolution of $\ZZ$
  over $\ZZ G$, and
  consider the commuting diagram of $\ZZ\gg$-modules
  $$
    \xymatrix{
    \ZZ &\ZZ\gg \ar[d]^{f_0}\ar[0,-1]_{\varepsilon}&&\ZZ\gg
    \ar[d]^{f_1}\ar[0,-2]_{g-1}& &\ZZ\gg 
    \ar[d]^{f_2}\ar[0,-2]_{\on{norm}} &&\ZZ\gg
    \ar[d]^{f_3}\ar[0,-2]_{g-1}  
    \\
    \ZZ&F_0\ar[l]&&F_1\ar[0,-2]_{\del_1}&&F_2\ar[0,-2]_{\del_2}
    &&{\phantom{,}F_3\ar[0,-2]_{\del_3}},    
    }
  $$
  where $$\on{norm}(1)=1+g+g^2+\dots+g^{n-1},$$ and
  in bar notation, the maps $f_i$ are given by
  $$
    f_0(1) = 1,\quad
    f_1(1) =[g],\quad
    f_2(1) = \sum_{k=0}^{n-1}[g^k|g],\quad\text{and}\quad
    f_3(1) = \sum_{k=0}^{n-1}[g|g^k|g].
  $$
  These are the first maps of an augmentation preserving map of chain
  complexes $f_\bullet$ from the periodic resolution to the standard
  resolution. On coinvariants, $f_\bullet$ induces an isomorphism of
  complexes of the form
  $$
    \ZZ\overset=\longleftarrow    \ZZ\overset0\longleftarrow
    \ZZ\overset n\longleftarrow    \ZZ\overset0\longleftarrow
    \ZZ\overset n\longleftarrow\cdots,
  $$
  where $1\in\ZZ\gg$ becomes $1\in\ZZ$.
\end{Pf}
\begin{Cor}\label{h-Cor}
  The order of the action $T^n$ on the $g$-twisted sector
  $$
    \bigoplus_{h\in C_g}L_{g,h}
  $$
  equals the order of $\alpha\at{\gg}$. In particular, it divides the
  order of $\alpha$ and that of $g$.
\end{Cor}
Following \cite{Conway:Norton}, we will write $h$ for the order of
$\alpha\at\gg$ and $n$ for the order of $g$, and set $N=nh$.
\begin{Cor}
  The map $f_g$ factors through
  \begin{eqnarray*}
    \HH&\to&\CC^\times\\
    \tau &\mapsto & e^{2\pi i\frac\tau {N}}.
  \end{eqnarray*}
\end{Cor}
\begin{Cor}
  There is a Fourier expansion
  $$
    f_g\in R(\hat C_g)\tensor\CC\ps{q^{\frac 1{N}}}
  $$
  of $f_g$, where $q=e^{2\pi i\tau}$.
  If $T^n$ acts by multiplication with $e^{2\pi i\frac sh}$ on
  $L_{g,1}$, then $f_g$ has the form
  $$
    f_g(\tau) = \sum_k q^{\frac kn+\frac s{N}}.
  $$
\end{Cor}
\begin{Pf}{}
  The first part follows immediately from the previous corollary. The
  second part is proved by a comparison of coefficients of the two
  sides of the equation
  $$
    f(q^{\frac1{N}}e^{2\pi i\frac1h}) = f(\tau+n) = f(\tau)e^{2\pi
      i\frac sh} = f(q^\frac1{N})e^{2\pi i\frac sh}.
  $$
\end{Pf}
In \cite{Dong:Li:Mason:modinv}
such a fraction $\frac s{N}$ turns up as the conformal weight of the
$g$-twisted sector.
If $g=1$, only integral powers of $q$ occur, and the central extension
of the centralizer is trivial:
$$\hat C_1=U(1)\times G.$$
\begin{Prop}\label{hth-root-Prop}
  The group $\Gamma_0(n)$ acts on $L_\gg$ by
  multiplication with $h^{th}$ roots of $1$. The group
  $\Gamma_0(N)$ acts trivially on $L_\gg$.
\end{Prop}
\begin{Pf}{}
  The first claim follows immediately from Lemma
  \ref{change-of-groups-Lem} and Lemma \ref{Gamma-0-Lem}. Let
  $\gamma\in\Gamma_0(N)$, and let $\map{p_h}{\TT^2}{\TT^2}$ be the
  degree $h$ map of the first circle. Then we have
  $$p_h\circ\gamma = \hat\gamma\circ p_h$$ with
  $\hat\gamma\in\Gamma_0(n)$. As in the proof of Lemma
  \ref{Gamma-0-Lem}, it follows that $\gamma_*=(\hat\gamma_*)^h$. 
  But we just saw that $\hat\gamma_*$ acts as an $h^{th}$ root of unity.
\end{Pf}
\noindent
\begin{Lem}\label{restriction-Lem}
  If $\left[\alpha\right]$ is an
  element of order $h$ in $H^3(B\gg,\mathbb R/\ZZ)$, then
  the restriction $\left[\alpha\at{B\langle 
    g^h\rangle}\right]$ has order $1$ in $H^3(B\langle g^h\rangle,\mathbb R/\ZZ)$.
\end{Lem}
\begin{Pf}{}
  There is a $\ZZ\langle g^h\rangle$-equivariant,
  augmentation-preserving map between the periodic resolutions,  
  which is the inclusion of $\ZZ\langle g^h\rangle$ in
  $\ZZ\gg$ in 
  even degrees and sends $1$ to $1+g+\dots+g^{h-1}$ in odd
  degrees. On coinvariants, this map becomes the identity in even
  degrees and multiplication with $h$ in odd degrees.
\end{Pf}

An alternative proof for this lemma is given as follows:
  Let $T^t$ denote the transpose of $T$.
  Note that we can use the degree $h$ map of the second circle to
  conjugate $(T^t)^{n}$ into $(T^t)^{\frac nh}$. It follows 
  that the action of $(T^t)^{\frac nh}$ on $L_{\langle g^h\rangle}$ is the
  $h^{th}$ power of the action of $(T^t)^n$ on $L_\gg$.
\begin{Cor}
  The group $\Gamma_0(\frac nh)$ acts trivially on $L_{\langle
    g^h\rangle} = L_\gg^{\tensor h}$.
\end{Cor}
\begin{Def}[{cf.\ \cite{Conway:Norton}}]
  Let $\Gamma_0(n|h)$ be the subgroup of $\on{SL}_2(\mathbb R)$
  consisting of matrices of the form
  $$
    \(\begin{array}[]{cc}
    a&\frac bh\\ 
    nc&d
    \end{array}\),
  $$
  where $a$, $b$, $c$ and $d$ are integers. 
\end{Def}
Note that $\Gamma_0(n|h)$ is conjugate to $\Gamma_0(\frac nh)$ via
conjugation with
$$
  \(\begin{array}[]{cc}
  h&0\\ 
  0&1
  \end{array}\)\inv\quad\text{or with}\quad
  \(\begin{array}[]{cc}
  1&0\\ 
  0&h
  \end{array}\).
$$
\begin{Rem}
The conjecture in \cite{Conway:Norton} suggests that 
$\Gamma_0(n|h)$ acts by $h^{th}$ roots of $1$ on $L_\gg$. While the
previous corollary does not tell us how to define this action when $h$
does not divide $b$, it shows
that in order for it to be consistent with the rest of the picture, it
needs to be an action by $h^{th}$ roots of unity.
\end{Rem}
\section{Abelian groups and level structures}
Let $E=\CC/\Lt$. Since any classifying map $\map f{\ZZ^2}G$ factors
through an abelian subgroup of $G$, it is often enough to study
principal bundles over $E$ with abelian structure group. These also
are important for the definition of Hecke operators.
\subsection{Isogenies}
Let $A$ be a finite abelian group.
Our first observation is that principal $A$-bundles with connected
total space are almost the same as isogenies with kernel $A$. More precisely, 
let $\map\phi {E'}E$ be an isogeny with kernel $A$. Then $A$ acts on
$E'$ by addition, making $\phi$ the structure map of a principal
$A$-bundle over $E$. 
\begin{Prop}\label{isogeny-Prop}
  Let $\map f{\ZZ^2}A$ be a group homomorphism. Then $f$ is surjective
  if and only if the principal $A$-bundle
  $\map\phi PE$ over $E$ classified by $\on{Bf}$ is isomorphic to
  an isogeny with kernel isomorphic to $A$, viewed as a principal
  $A$-bundle in the way just described. This in turn is the case if
  and only if $P$ is connected. 
\end{Prop}
%
\begin{Pf}{}
  Let $\map f{\ZZ^2}A$ be surjective and let $n$ be the order of $A$.
  Write $\Lambda$ for $\ker(f)$. Then $\Lambda$ is an index $n$
  subgroup of $\langle\tau,1\rangle$, and therefore a rank two
  sublattice. The short exact sequence
  $$0\to\Lambda\to\ZZ^2\to A\to 0$$
  gives rise to a fibre sequence of classifying spaces
  $$\on{B\Lambda}\longrightarrow\BZZ
  \stackrel{\on{Bf}}{\longrightarrow}BA,$$  
  whose first map is the principal $A$-bundle classified by
  $\on{Bf}$. 
  It is constructed as the quotient map
  \begin{equation}
    \label{isogeny-Eqn}
    \on{B\Lambda}:=\on{E\ZZ^2}/\Lambda =
    \CC/\Lambda\to\CC/\langle\tau,1\rangle = \BZZ,    
  \end{equation}
  with the $A$ action induced by the action of
  $\langle\tau,1\rangle\cong\ZxZ$ on
  $\CC\cong\on{E\ZZ^2}$. 
  Hence it is an isogeny whose kernel is isomorphic to $A$.
  Moreover, any such isogeny into $E$ is of the form
  (\ref{isogeny-Eqn}) for some sublattice $\Lambda\sub\Lt$, and
  therefore classified by the surjective map 
  $$\ZxZ\cong\Lt\twoheadrightarrow\Lt/\Lambda.$$
  Let now $P$ be connected. We will see in Lemma \ref{connected-Lem}
  that this implies that $P$ can be viewed as a complex elliptic curve and
  $\phi$ as an isogeny.
  More precisely, $P$ is the quotient of $\CC$ by an index $n$ subgroup
  $\Lambda\sub\Lt$. Thus
  $\phi$ is an isogeny with kernel $\Lt/\Lambda$, and Proposition \ref{GHisom-Prop}
  implies that $A$ is isomorphic to $\Lt/\Lambda$.  
\end{Pf}
\begin{Rem}\label{isogenies-neq-principal-Rem}
  The (base-preserving) automorphisms of a principal $A$-bundle over $E$
  are given by 
  multiplication with elements of $A$:
  $$
    \Aut(P_{a_1,a_2})\cong C_A(a_1,a_2) = A.
  $$
  Isogenies, at the other hand, come equipped
  with a choice of basepoint over zero and therefore have no
  automorphisms over $E$.
\end{Rem}

If $f$ is not surjective, we can write $f=i\circ f'$, where
$\map {f'}{\ZZ^2} B$ is surjective, and $\map i{B}A$ is injective.
Then the $A$-bundle classified by $\on{Bf}$ is a collection of
isogenies into $E$ with kernel $B$.
\subsection{Duals}
Instead of considering isogenies into $E$, it is often more natural to
consider the dual picture, which involves isogenies with source $E$ and the
corresponding level structures on $E$. 
Let $\hat A:=\Hom(A,\bbS^1)$ be the Pontrjagin dual of $A$. Since
$A$ is finite and abelian, we have $\hat{\hat{A}}=A$. 
We further know $$\widehat\ZZ=\bbS^1 \text{ and } \widehat{\bbS^1}=\ZZ$$
Let $i_{\langle\omega_1,\omega_2\rangle}$ denote the isomorphism from
$E$ to $\bbS^1\times\bbS^1$ corresponding to the choice of basis
$(\omega_1,\omega_2)$. 
While we had identified $$B\ZZ^2=\bbS^1\times\bbS^1$$ with $E$ in the
non-standard way, using $i_{\langle-1,\tau\rangle}$, we will follow
the convention to identify
$$\widehat{\ZZ^2}\cong\bbS^1\times\bbS^1$$ with $E$ in the standard way,
using $i_{\langle\tau,1\rangle}$. 
\begin{Prop}[{compare \cite[(1.4.4)]{Ginzburg:Kapranov:Vasserot}}]
These identifications yield an isomorphism
$$\Prin A{E}\cong\Hom(\hat A,E),$$
where $\Prin AE$ denotes the set of isomorphism classes of principal
$A$-bundles over $E$. 
Let $\map\phi{E'}E$ be an isogeny.
Then pullback along $\phi\dual$ on the left-hand side 
becomes composition with $\phi$ on the right-hand side.
\end{Prop}
\begin{Pf}{}
  The isomorphism is given by
  $$\Prin AE\cong\Hom(\ZxZ,A)\cong\Hom(\hat A,\widehat{\ZZ^2})
     \cong \Hom(\hat A,E)
  .$$
  Here the first isomorphism uses 
  $i_{\langle -1,\tau\rangle}$
  to identify $\bbS^1\times
  \bbS^1$ with $E$, while the last one uses
    $i_\Lt$.  
  Let $\map\varphi{\ZZ^2}{\ZZ^2}$ be as in (\ref{phi-Eqn}) (on page
  \pageref{phi-Eqn}), let 
  $\varphi\dual$ denote its 
  pseudo-inverse, and write $\phi=B\varphi$. Then $\phi$ is viewed as a
  map from $E'$ 
  to $E$ via the identifications $i_{\langle -1,\tau'\rangle}$
  and $i_{\langle -1,\tau\rangle}$. 
  The first isomorphism sends 
  $(\phi\dual)^*$ to $(\varphi\dual)^*$. The second isomorphism maps
  $(\varphi\dual)^*$ to $\hat{\phi\dual}_*$.
  We have a commutative diagram
  \begin{equation}
    \label{commutative-diagram-dual-Eqn}
    \xymatrix{
    {E'{\phantom{}}}\ar[0,2]^{i_\Ltp{\phantom{xxx}}}\ar[d]_{\phi} && 
    {{\phantom{}}\bbS^1\times\bbS^1}\ar[0,2]^{{\phantom{xx}}\cong}
    \ar[d]^{B(\varphi\dual)^t}&& {\hat{\ZZ^2}}   
    \ar[d]^{\hat{\varphi\dual}}
    \\
    {E}\ar[0,2]^{i_\Lt{\phantom{xxx}}}& &
    {\bbS^1\times\bbS^1}\ar[0,2]^{{\phantom {xx}}\cong} &&{\hat{\ZZ^2}},
    }
  \end{equation}
  where $(\varphi\dual)^t$ is
  the transpose of $\varphi\dual$.
  This completes the proof.
\end{Pf}
\begin{Lem}\label{pullback-Lem}
  Let $f\in\Hom(\hat A,E)$, and 
  let $\xi_f$ be the principal
  $A$-bundle over $E$ corresponding to $f$.
  Let $g$ be an automorphism of $A$. Then the pullback
  $g^*(\xi_f)$ of $\xi_f$ along $g$ is the principal $A$-bundle
  corresponding to $f\circ\widehat{g^{-1}}\in\Hom(\hat A,E)$,
  $$g^*(\xi_f)\cong\xi_{f\circ\widehat{g^{-1}}}.$$
\end{Lem}
\begin{Pf}{}
  Let $\map{\hat f}{\ZxZ}{A}$ be the dual of $f$. Then $\xi_f$ is
  classified by $\on{B\hat f}$, and $\xi_{f\circ\widehat{g^{-1}}}$ is
  classified by
  $$\on{B(\widehat{f\circ\hat{g^{-1}}})}=\on{B(g^{-1}\circ 
    \widehat{f})}=\on{Bg^{-1}}\circ\on{B\hat f}.$$  
  The proposition now follows from Corollary \ref{Balpha-Cor}.
\end{Pf}

A homomorphism $\map f\ZxZ A$ is surjective, if and only if its dual
homomorphism $\map{\hat f}{\hat A}E$ is injective.
\begin{Prop}
  Let $f\in \Inj(\hat A,E)$, and let $\xi_f$ be the principal
  $A$-bundle over $E$ classified by $B\hat f$.
  Further, let $\map{\Psi_f} EE'$ be the
  isogeny with kernel $f$.
  Then $\xi_f$ is isomorphic to the principal bundle defined by
  the dual isogeny $\map{\Psi\check{{}_f}}{E'}{E}$. 
\end{Prop}
\begin{Pf}{}  
  Consider the short exact sequence
  $$0\to\hat A\stackrel{f}\longrightarrow E\stackrel{\Psi_f}\longrightarrow 
  E'\to 0.$$ 
  Its Pontrjagin 
  dual is a short exact sequence of the form
  $$0\to \ker(\hat f)\longrightarrow\ZxZ\stackrel{\hat f}{\longrightarrow} A\to 0.$$
  Write $\psi$ for the inclusion of
  $\ker(\hat f)$ in $\ZxZ$.
  Then $\psi$ is such that under the isomorphisms $i_\Lt$ and $i_\Ltp$, we have
  $\Psi_f=\hat\psi$.   
  It follows from the proof of Proposition \ref{isogeny-Prop} that
  under the isomorphisms $i_{\langle-1,\tau\rangle}$ and $i_{\langle
    -1,\tau'\rangle}$, we have $\xi_f=B\psi$.  
    Finally, we consider the commutative diagram
  (\ref{commutative-diagram-dual-Eqn})
  with $\xi_f\check{}$ in the role of $\phi$ and $\psi$ in the role of
  $\varphi\dual$ to obtain $\xi_f\check{}=\Psi_f$.
\end{Pf}
\begin{Cor}\label{Inj-Cor}
  There is a one-to-one correspondence
  $$\Inj(\hat A,E)\longleftrightarrow\{\text{principal $A$-bundles
    over $E$ with connected total space} \}\medspace\slash\medspace\cong$$
  sending $f$ to $\Psi\check{{}_f}$. 
\end{Cor}
The decomposition
$$
  \Hom(\hat A,E)\cong\coprod_{B\sub A}\Inj(\hat B,E)
$$
is invariant under the action of $\SL$ on $A\times A$. Hence it
induces a decomposition of the coarse moduli space
$$
  X_A = \coprod_{B\sub A} \XinjB,
$$
where $\XinjB$ is the moduli space of elliptic curves with level
$\hat B$-structures, i.e., of isogenies $\map\phi EE'$ together with a
choice of isomorphism from their kernel to
$\hat B$. 
Note that $$\Inj(\hat A,E)$$ is
invariant under the action of $\Aut(A)$ 
on $\HomAE$, and that
$$\Inj(\hat A,E)/\Aut(A)$$
corresponds to the set of all isogenies out of $E$ which allow a
choice of isomorphism from 
$\hat A$ to their kernels. Since no choice of this isomorphism is
specified, and $\hat A$ is non-canonically
isomorphic to $A$, this condition is equivalent to allowing an
isomorphism from $A$ to their kernels.
\begin{Def}
  We will write $\MconnA$ for the component of $\mM_A$ parametrizing
  bundles with connected total space and  
  $$
    \MzeroA := \mM_A^{\on{conn}}\mmod \Aut(A)
  $$
  for the quotient stack of $\MconnA$ 
  by the action of $\Aut(A)$ on
  the fibers. We will write $X_0(A)$ for the corresponding coarse
  moduli space. 
\end{Def}
Note that by Remark \ref{isogenies-neq-principal-Rem} the stack
$\MzeroA$ is not the same as 
the moduli-stack $\mM_0(A)$ of isogenies allowing an isomorphism from
$A$ to their kernel.
However, the forgetful 
map 
$$
  \mM_0(A)\longrightarrow\MzeroA
$$
induces an isomorphism of coarse moduli spaces.
\begin{Def}
  The Fricke involution $W_A$ is an involution of $\mM_0(A)$.
  It sends the isogeny $\phi$ to its dual isogeny
  $$
    \map{\phi\dual}{E'}{E}.
  $$
  If the kernel of $\phi$ is isomorphic to $\hat A$, then the kernel of
  $\phi\dual$ is (also) isomorphic to $A$.   
\end{Def}
\subsection{Cyclic groups}
Let now $A$ be the cyclic group $\ZZn$.
An injective map from $\ZZn$ to $E$ is the same as a choice of point
of exact order $n$ on $E$. Such a choice is called a $\Gamma_1(n)$-level
structure on $E$. Therefore, the space of principal $\ZZn$ bundles
on $E$ is non-canonically isomorphic to
$$E[n]\cong\Hom({\ZZn},E)\cong
\coprod_{d|n}\on{Level}_{\Gamma_1(d)}(E),$$
and the moduli space $X_\ZZn$ gets (non-canonically) identified with
$$
  X_\ZZn \cong \coprod_{d|n} X_1(d).
$$
Note that
$$
  \Inj(\ZZn,E)/\Aut(\ZZn) \cong \on{Level}_{\Gamma_0(n)}(E).
$$
Therefore, we get a (canonical) isomorphism
$$
  X_\ZZn^\inj /\Aut(\ZZn) \cong X_0(n).
$$
The Fricke involution $W_n:=W_A$ takes a cyclic subgroup $C\sub E$
of 
order $n$ to the cyclic subgroup
$$
  E[n]/C\sub E/C.
$$
Let $g\in\ZZ/n\ZZ$ be a generator. Then the triple $(1,g;\tau)$ corresponds
to the cyclic subgroup $\left\langle\frac 1n\right\rangle$ of
$\CC/\left\langle\tau,1\right\rangle$. Therefore $W_n$ sends it to the cyclic
subgroup of $\CC/\left\langle\tau,\frac 1n\right\rangle$ generated by
$\frac\tau n$. This again corresponds to the subgroup generated by
$\frac 1n$ in 
$$\CC/\left\langle\frac{-1}{n\tau},1\right\rangle.$$ 
Thus
$$
  W_n\(1,g;\tau\) = \(1,g;\frac{-1}{n\tau}\).
$$
\section{Symmetric groups and coverings}
The symmetric groups play a central role in the definition of power
operations and will turn out to be important in the theory of
replicability.  
Recall the one-to-one
correspondence 
\begin{eqnarray*}
  \{\text{Isom.\ classes of principal $\Sn$-bundles over $X$}\}&
  \longleftrightarrow&
  \{\text{Isom.\ classes of $n$-fold covers of $X$}\} \\
  P &\mapsto& P\times_\Sn \underline{n},
\end{eqnarray*}
where $\underline n$ is a fixed set with $n$ elements.
\begin{Lem}\label{connected-Lem}
Let $\map\pi PE$ be an $n$-fold covering of a complex elliptic
curve $E$ with connected
total space. Then $P$ can be viewed as a complex elliptic curve and
$\pi$ as an isogeny.
\end{Lem}
\begin{Pf}{}
By the theory of covering spaces, the universal cover 
of $E$,
$$\CC\to\CC/\Lt\cong E,$$
factors through $\pi$.
More precisely, if we choose a basepoint of $P$ over zero, this
identifies $P$ with
the quotient of $\CC$ by an index $n$ subgroup
$\Lambda\sub\Lt$. Thus
$\pi$ becomes an isogeny.  
\end{Pf}

It follows that every $\Sn$-bundle with connected total space is
induced from a principal bundle with abelian structure group:
Let $A$ be an abelian group of order $n$, and
pick a bijection between $A$ and $\underline n$. Then the action
of $A$ on itself by multiplication defines a map 
$$
  i\negmedspace :A\hookrightarrow\Sn.
$$
Note that a different choice of bijection leads to a different map which is
conjugate to $i$ by an element of $\Sn$. 
In particular, this applies to precomposition of $i$ with 
any automorphism of $A$.
\begin{Lem}
Let $$\map\phi Y{E}$$ be an isogeny together with an isomorphism
from $\ker(\phi)$ to $A$.
Consider the $n$-fold covering map underlying $\phi$, and let $\xi_Y$ be
the corresponding principal $\Sn$-bundle. Now consider $\phi$ as
principal $A$-bundle. Then the isomorphism class of the
$\Sn$-bundle $\phi[\Sn]$ associated to
$\phi$ via $i$ is equal to that of $\xi_Y$. 
\end{Lem}
Note that (in particular) the isomorphism class of $\phi[\Sn]$ is
independent of the choice of $i$ and of the choice of isomorphism
$\ker(\phi)\cong A$.

\bigskip
\begin{Pf}{}
  By Lemma \ref{associated-fiber-Lem} we have 
  $$\phi[\Sn]\times_\Sn\underline n\cong\phi\times_A\underline n,$$
  which is the $n$-fold covering underlying $\phi$. 
  Thus $\xi_Y$ and $\phi[\Sn]$ correspond to the same covering, which
  proves the claim.
\end{Pf}
%
Take an $n$-fold covering $P$ of $E$ with connected total space, and view
it as a principal $A$-bundle. Note that any deck transformation can be
realized as $(-)+a$ with $a\in A$. Hence the automorphisms of $P$ as
an $n$-fold covering are the same as the automorphisms of $P$ as
principal $A$-bundle. For future reference, we note that their number
is $n=|A|$.

We now direct our attention to coverings whose total spaces are not
connected. 
Recall that
if $P$ and $P'$ are $\Sn$- and $\Sm$-principal bundles with associated
covering spaces $Y$ and $Y'$ respectively, then the associated
covering space of 
$$
  (P\times P')\times_{\Sn\times\Sm}\Snm
$$
is the disjoint union $Y\cup Y'$.
Let now $Y$ be the disjoint union of connected components
$$
{Y = Y_1}\cup\cdots \cup{ Y_k}.
$$  
Then the $Y_i$ can be viewed as having abelian structure groups
$A_1,\dots,A_k$ of orders 
$n_1,\dots,n_k$ 
with $n=\sum n_i$. 
Each of the $Y_i$ has an associated $\Sigma_{n_i}$-principal bundle 
$$
  P_i = Y_i\times_{A_i}\Sigma_{n_i},
$$
and $P$ is the bundle
$$
  P = \(\prod_iP_i\)\times_{\Sigma_1\times\dots\times\Sigma_k}\Sn.
$$
Note that each $A_i$ is a quotient of $\ZZ^2$, making
the disjoint union
$$
  A_1\cup\dots\cup A_k\cong\underline{n},
$$
into an $\ZZ^2$-set,
whose orbits are the $A_i$. The $\ZZ^2$-action is transitive if and
only if $Y$ is connected.
We obtain a map 
$$
  \ZZ^2\to\Sn,
$$
well defined up to conjugation in $\Sn$, and this is the pair of
commuting elements which classifies $P$.
More precisely, we have the following:
\begin{Prop}\label{Sn-decomp-Prop}
The decomposition of covers into connected components defines an equivalence 
$$
  \isomap{e}{\mM_\Sn}{
    \coprod_{n=\sum|A|N_A}\prod_{[A]}\MzeroA\wr\Sigma_{N_A}}
$$  
where the product as well as the sum under the coproduct is over all
isomorphism classes of 
finite abelian groups, and all products are taken over the
target maps $\widetilde\mM_A\to\mM_1$.
\end{Prop}
\begin{Exa}
  Let $n=9$, and consider the decomposition $9 = |\{1\}| +
  |\ZZ/4\ZZ|\cdot 2$. The corresponding summand 
  $$\mM_1\times_{\mM_1}\widetilde \mM_{\ZZ/4\ZZ}\wr\Sigma_2$$
  parametrizes ninefold covers of elliptic curves that decompose into three
  connected components, two of which are fourfold covers corresponding
  to isogenies with kernel $\ZZ/4\ZZ$ and one of which is a onefold cover.  
\end{Exa}
\begin{Pf}{of Proposition \ref{Sn-decomp-Prop}}
  We start by considering bundles over a torus. 
  Fix a partition $$n=\sum_T|T|N_T,$$ where the $T$ are elements of a
  fixed system of representatives of the set of isomorphism classes of
  finite transitive $\ZZ^2$ sets.
  Such decompositions classify conjugacy classes $[\sigma,\rho]_\Sn$
  of pairs of commuting elements in $\Sn$, and hence they classify 
  isomorphism classes of principal $\Sn$-bundles over $\TT^2$. 
  For each $T$, let 
  $$A_T:= \ZZ^2/\Stab_{\ZZ^2}(T)$$ 
  denote the corresponding abelian group.
  Note that the isomorphism class of $T$ is uniquely determined by
  the surjective map  $\map{f_T}{\ZZ^2}{A}$, that for any automorphism
  $g$ of $A$, the map $g\circ f_T$ determines the same isomorphism
  class $[T]$, and that we have
  $$N_A=\sum_{A_T= A}N_T.$$ 
  Now the $n$-fold cover of $\TT^2$ corresponding to $[\sigma,\rho]_\Sn$
  decomposes into connected
  components,  
  $N_T$ of which are isomorphic to the principal $A_T$-bundle with
  connected total space that is classified by $f_T$. 
  The automorphism group of this
  $n$-fold cover over $\mathbb T^2$ is isomorphic to 
  $$
    C_\Sn(\sigma,\rho) = \prod_TA_T^{N_T}\wr\Sigma_{N_T}.
  $$
  For each $[A]$, pick, once and for all, an ordering of the classes
  $[T]$ with $A_T\cong A$. Then, for each pair $[\sigma,\rho]_\Sn$ and for each
  orbit-type $[T]$ occuring in its orbit decomposition, pick an
  ordering of the orbits isomorphic to $T$. 
  In this way, we assign to $[\sigma,\rho]_\Sn$ a point  
  $$
    e\st \in \prod_A \(\on{Surj}(\ZZ^2,A)/\Aut(A)\)^{N_A}.
  $$
  Our choices also determine an isomorphism  
  $$\Aut(e\st)\cong C_\Sn\st,$$
  which we use to continue $e$ to automorphisms of $[\sigma, \rho]$.
  Let now $\gamma\in\SL$. Then $\gamma$ acts from the right on $\st$
  and on $\on{Surj}(\ZZ^2,A)$, for each $A$. The transitive
  $\ZZ^2$-set $T\cdot\gamma$ classified by $\gamma\circ f_T$ might be different
  from $T$, but we have
  $$
    A_{T\cdot\gamma} = A_T.
  $$
  The action of $\gamma$ on $\st$ is compatible with the action of
  $\st$ on the connected components of the corresponding
  bundle. However, we might have numbered the components 
  belonging to
  the same group $A$ differently for $\st\cdot\gamma$ than we did for
  $\st$, with the result that they differ by an element
  $\sigma_\gamma\in \Sigma_{N_A}$. 
  To extend the map $e$ to the translation groupoid
  $$
    \(\ol\mC(\Sn)\times\HH\)\rtimes\SL,
  $$
  we set 
  $$
    e((\sigma,\rho;\tau),\gamma) := 
    \prod_A((\st_A;\tau),(\gamma\wr\sigma_{\gamma})).
  $$ 
  Here $\st_A$ stands for the collection of the $\st$-orbits that
  belong to $A$.

  By construction, $e$ is fully faithful and essentially
  surjective. Different choices of orderings would result in a map
  $e'$, which 
  would differ from $e$ by a natural transformation on its target.
  Hence, as a map of 
  orbifolds, $e$ is independent of these choices. 
\end{Pf}
\begin{Def}\label{Hecke-correspondence-Def}
  We will write $\widetilde{\mH}_n$ for the moduli stack of all $n$-fold covers
  with connected fibre,
  $$
    \widetilde{\mH}_n = \coprod_{|A|=n}\MzeroA.
  $$ 
  Here the disjoint union runs over all isomorphism classes of abelian
  groups $A$ of order $n$.  
\end{Def}

Equivalently, $\widetilde{\mH}_n$ is the substack of $\mM_\Sn$ defined as
$$
  \widetilde\mH_n:=\{f:\ZZ^2\to\Sn \mid \text{$f$ makes $\underline
    n$ into a 
  transitive $\ZZ^2$-set}\}\times\HH\medspace\mmod\medspace (\Sn\times\SL).
$$

The following lemma provides an alternative proof for the fact that the
automorphism group of an $n$-fold covering with connected total space
is isomorphic to the kernel $A$ of the corresponding isogeny.
\begin{Lem}
  Let $(\sigma,\rho)$ be a pair of commuting elements of $\Sn$ which
  acts transitively on $\{1,\dots,n\}$. Write
  $A:=\langle\sigma,\rho\rangle$. Then $A$ has order $n$, the
  principal $\Sn$-bundle classified by $\st$ corresponds to an isogeny
  with kernel $A$, and we have 
  $$C_\Sn(\sigma,\rho) = \langle\sigma,\rho\rangle = A.$$
\end{Lem}
\begin{Pf}{}
  If $\st$ acts transitively on $\{1,\dots,n\}$, then the
  corresponding $n$-fold covering has connected total space and hence
  is an isogeny. Write $A'$ for its kernel, then the bundle $P_\st$ is
  induced by a principal $A'$-bundle. Hence the classifying map of
  $P_\st$ factors through $A'$, yielding an isomorphism to $A$:
  $$\ZxZ\to A'\stackrel{\cong}\longrightarrow
  \langle\sigma,\rho\rangle\to\Sn.$$ 
  To determine the centralizer of
  $\langle\sigma,\rho\rangle$, 
  let $\pi\in C_\Sn(\sigma,\rho)$. 
  Because of the transitivity of the action of $A$, we can write $\pi(1)$ in
  the form $\pi(1)=\sigma^a\rho^b(1)$. Further, for $1\leq x\leq n$,
  we have $c$ and $d$ such that $x=\sigma^c\rho^d(1)$. It follows that 
  $$
    \pi(x) = \pi(\sigma^c\rho^d(1)) = \sigma^c\rho^d(\pi(1)) = 
    \sigma^c\rho^d(\sigma^a\rho^b(1)) =
    \sigma^a\rho^b(\sigma^c\rho^d(1)) = \sigma^a\rho^b(x).    
  $$
  Hence $\pi= \sigma^a\rho^b$.
\end{Pf}
\section{\HKR character theory}
Let $\mA(G)$ be the category having as objects the abelian subgroups of
$G$ and with morphisms from $B$ to $A$ being the $G$-equivariant maps
from $G/B$ to $G/A$. Then $\mA(G)$ is a full subcategory of the
standard orbit category. 
The results of the previous
sections can be summarized as follows: 
\begin{eqnarray}\label{Artin-Eqn}
  \mM_G & \simeq & \Colim_{\mA(G)} \mM_A.
\end{eqnarray}
Here $\Colim{}$ stands for a weak $2$-colimit. Roughly speaking, that
means that we take the disjoint union of all the $\mM_A$ and instead
of quotienting by the equivalence relation defined by
maps in $\mA(G)$, as we would for a classical colimit, we add in extra
isomorphisms.

In fact, one can
replace $\mA$ by the full subcategory whose objects are abelian groups
which are generated by two elements.
The equivalence \eqref{Artin-Eqn} is an analogue of the generalized
Artin theorem 
\cite[Thm. A]{Hopkins:Kuhn:Ravenel}. 
\subsection{Products} 
Let $G$ and $H$ be finite groups. Then we have an isomorphism of
moduli stacks
$$
  \mM_{G\times H}\stackrel{\cong}\longrightarrow \mM_G\times_{\mM_1}\mM_H,
$$
which pulls back $\mL^\alpha\tensor \mL^\beta$ to $\mL^{\alpha\times\beta}$.
This yields external products
$$
  \map\boxtimes{
  \Ell^\alpha_G(\pt)\tensor  \Ell^\beta_H(\pt)}{
  \Ell^{\alpha\times\beta}_{G\times H}(\pt)},
$$ 
$$
  (f_G\boxtimes f_H)(g_1,h_1,g_2,h_2;\tau) = f_G(g_1,g_2)\cdot f_H(h_1,h_2).
$$ 
In the case that $G=H$, one composes with pull-back along the diagonal
of $G\times G$ to obtain internal products
$$
  \map\tensor{
  \Ell^\alpha_G(\pt)\tensor  \Ell^\beta_G(\pt)}{
  \Ell^{\alpha+\beta}_{G}(\pt)}.
$$ 
\subsection{Change of groups}
Let $\map a HG$ be a map of finite groups.
We write 
$$\map{\overline{\mathcal C}(a)}{\overline{\mathcal
  C}(H)}{\overline{\mathcal C}(G)}$$  
for the map of moduli sets induced by $a$, and note that
there is a canonical isomorphism of line bundles
$$
  \overline\mC(a)^*(L^\alpha)\cong L^{a^*(\alpha)}.
$$
For a section $f$ of $L^\alpha$ over $\overline{\mathcal C}_G$, we
define the section $\res{}a(f)$ of
$L^{a^*(\alpha)}$
over $\overline\mC(H)$ by  
$$
  \res{}a(f) := f\circ\overline{\mathcal C}(a).
$$
We recall the measure $\mu$ on $\overline{\mathcal C}(G)$ defined in
\cite[(2.1)]{Freed:Quinn}: 
$$
  \mu([P_{g,h}]) = \frac1{|\Aut(P_{g,h})|} = \GG{C_G(g,h)}.
$$
If $f$ is a section of $L^{a^*(\alpha)}$ over
$\overline{\mathcal C}(H)$,
we define a section $\ind a{}(f)$ of $L^\alpha$
over $\overline\mC(G)$ 
by 
$$
  \ind a{}(f)\cdot\mu = \int_{fiber} f\cdot d\mu = \sum_{fiber}
  f([h_1,h_1]_H)\GG{C_H(h_1,h_2)}, 
$$
where at the point $[g_1,g_2]_G$, the sum is over all $H$-conjugacy
classes $[h_1,h_2]_H$ which 
satisfy $[h_1,h_2]_G = [g_1,g_2]_G$.
We have
\begin{eqnarray*}
  \ind a{}(f)([g_1,g_2]_G) & = & 
|C_G(g_1,g_2)|\cdot\sum_{[h_1,h_2]_H\in fiber}
\GG{C_H(h_1,h_2)}\cdot f([h_1,h_1]_H)
\\
&=&|C_G(g_1,g_2)|\cdot\sum_{a(h_1,h_2)\sim_G(g_1,g_2)}\GG{H}\cdot f([h_1,h_1]_H).
\end{eqnarray*}
Two special cases are important: in the case that $a$ is an inclusion
of groups, we obtain the formula of \cite[Thm. D]{Hopkins:Kuhn:Ravenel}, 
$$
  \ind HG(f)
  =\GG H\sum_{\stackrel{s\in G\mid}{ s\inv(g_1,g_2)s\in H^2}}
  f(s\inv g_1 s, s\inv g_2 s).
$$
In the case that $a$ is the unique map from $G$ to the trivial group,
we get 
$$
  \ind G1(f)(1,1)= \GG G\sum_{gh=hg}f(g,h),
$$
where the sum is over all pairs of commuting elements of $G$.
\begin{Def}
  We define the inner product on $\Ell^0_G(\pt)$ as the composite of the
  (internal) product with $\ind G1$,
  $$
    \langle-,-\rangle_G:
    \Ell^0_G(\pt)\tensor\Ell^0_G(\pt)\to\Ell^0_G(\pt)
    \longrightarrow\Ell^0_1(\pt).
  $$
\end{Def}
We have
$$
  \langle f_1,f_2\rangle_G = \GG
  G\sum_{gh=hg}f_1(g,h)\cdot f_2(g,h).
$$
All of these constructions extend to $\SL$-equivariant maps from $\HH$
to $\bigoplus L_{g,h}$.
\section{Hecke operators and power operations}
\subsection{The Hecke correspondence}\label{Hecke-correspondence-Sec}
Let $\widetilde{\mH}_{n,G}\sub \mM_{G\times\Sn}$ be defined by
$$
  \widetilde{\mH}_{n,G}:= \mM_G\times_{\mM_1}\widetilde{\mH}_n,
$$ 
where $\widetilde{\mH}_n$ is as in Definition \ref{Hecke-correspondence-Def}.
Explicitly, we have 
\begin{eqnarray*}
  \widetilde{\mH}_{n,G} 
  & = & {\left\{\raisebox{3.5ex}{\xymatrix@=2ex{P'\ar[r]\ar[d] & 
          P\ar[d] & P' = \phi^*(P) \\
          E'\ar[r]^\phi&E&\deg(\phi)=n}}  
    \right\}\mmod\cong}, 
\end{eqnarray*}
where $P$ and $P'$ are principal $G$-bundles and $\phi$ is an $n$-fold
covering with connected total space.
Consider the diagram
\begin{equation}\label{Hecke-correspondence-Eqn}
  \xymatrix{
  &\widetilde{\mH}_{n,G}\ar[ld]_s\ar[rd]^t &\\
 \mM_G&&\mM_G,
  }
\end{equation}
where
the map $s$ returns $P'\to E'$, and $t$ returns $P\to E$.
Note $s$ is only defined canonically up to canonical natural
isomorphism in $\mM_1$. 
The map $t$ is induced by the map of groups
$$\map{\id\times p}{G\times\Sn}G,$$
where $p$ is the unique map from $\Sn$ to the trivial group.

Note that this 
is not 
the usual Hecke correspondence; like
$\MzeroA$, the stack $\widetilde\mH_n$ has more automorphisms than $\mH_n$ but
induces the correct coarse moduli space.
Therefore,
  in the case that $G$ is the trivial group, the diagram of maps
  between coarse moduli spaces induced by
  (\ref{Hecke-correspondence-Eqn})
  is the classical Hecke correspondence.
\begin{Def}
  We define
  the trace map $t_!$ as the restriction of
  $\ind{\id\times p}{}$ to $\widetilde{\mH}_{n,G}$.
\end{Def}
Proposition \ref{Ln-Prop} yields an isomorphism 
$$
  s^*\mL^\alpha\cong t^*(\mL^\alpha)^{\tensor n}
$$
over $\widetilde{\mH}_{n,G}$.
\begin{Def}\label{Hecke-operator-Def}
The $n^{th}$ Hecke operator 
$$
  \map{T_n}{\Gamma(\mL^\alpha)}{\Gamma(\mL^{n\alpha})}
$$
is defined as the composition
$$
  T_n = t_!\circ s^*.
$$
\end{Def}
\begin{Prop}\label{Hecke-Prop}
  Let $f$ be a section of $\mL^\alpha$. Then $T_n(f)$ is the
  following section of $(\mL^\alpha)^{\tensor n}$:
  $$
    T_n(f)(g,h;\tau) = \frac 1n\sum_{\stackrel{ad=n}{0\leq b<d}}
    \(\phi_{a,b,d}\)_*f\(g^d,g^{-b}h^a;\frac{a\tau+b}{d}\),
  $$
  where 
  $$\map{\phi_{a,b,d}}{\CC/\Ltp}{\CC/\Lt}
  $$ 
  sends $\tau'=(a\tau+b)/d$ to $a\tau+b$ and $1$ to $d$, and
  $$
    \map{(\phi_{a,b,d})_*}{ L^\alpha_{g^d,g^{-b}h^a}}{(L_{g,h}^\alpha)^{\tensor n}}
  $$
  is the induced isomorphism of line bundles of Proposition \ref{Ln-Prop}.  
\end{Prop}
\begin{Pf}{}
Let $P$ be a principal $G$-bundle over $E\cong\CC/\Lt$, let
$E'\cong\CC/\langle\tau',1\rangle$, where $\tau'=\frac{a\tau+b}{d}$ and let
$\phi=\phi_{a,b,d}$.
Then $\phi=B\varphi$, where 
\begin{eqnarray*}
\varphi\negmedspace : {\ZZ^2}&\to& {\ZZ^2}\\
(1,0)&\mapsto& (d,0)\\
(0,1)&\mapsto& (-b,a).  
\end{eqnarray*}
Let $P$ be the principal bundle over $E$ classified by the pair $(g,h)$.
Then the principal bundle $\Phi^*P$ over $E'$ is classified by the
pair $(g^d,g^{-b}h^a)$.  
\end{Pf}
\begin{Exa}
On the untwisted sector $g=1$, the formula of Proposition
\ref{Hecke-Prop} specializes 
to the formula for the twisted Hecke operators of classical Moonshine
in (\ref{twisted-Hecke-Eqn}).   
\end{Exa}
\begin{Exa}
  Let $G=\ZZl$. Then  
  $$
    (e,f)\in(\ZZl)^2
  $$
  corresponds to the point $\frac{e\tau+f}{l}$ of $E[l]$.
  Let $\omega_1$, $\omega_2$ and $\phi$ be as above.
  We have
  $$
    \frac{e\tau+f}{l} = de\frac{\omega_1}{l} - be\frac{\omega_2}{l} +
    af\frac{\omega_2}{l}. 
  $$
  Hence the image of $\frac{e\tau+f}{l}$ under $\phi$ is the point 
  of $E'[l]$ corresponding to $(de,af-be)$.
\end{Exa}
The action of the
twisted Hecke operators on the $q$-expansions of the untwisted sectors is
computed as follows:
\begin{eqnarray*}
  \tHn f(\tau) & = & \frac1n \sum_{\stackrel{ad=n}{0\leq
    b<d}} f^{(a)}\(\frac{a\tau+b}d\)\\ 
  & = & \frac1n\sum_{\stackrel{ad=n}{0\leq b<d}}
    \sum_{m\in\ZZ}c^{(a)}(m)q^{\frac{am}{d}}\zeta_d^{bm}\\
  & = &  \frac1n\sum_{ad=n}
    \sum_{m\in d\ZZ}d\cdot c^{(a)}(m)q^{\frac{am}{d}} \\
  & = & \sum_{m\in\ZZ}\sum_{a|(m,n)}\frac1a\cdot c^{(a)}\(\frac{nm}{a^2}\)q^m.
\end{eqnarray*}

\subsection{Power operations in K-theory}
In \cite{Atiyah:poweroperations}, Atiyah developed the theory of
cohomology operations in equivariant $K$-theory (also known as {\em
  equivariant power operations}). In this section, we
briefly recall the special cases of Atiyah's definitions that will
serve as a model for our definitions in $\Ell_G^\alpha(\pt)$.
Recall that the coefficient ring of equivariant $K$-theory is the
representation ring $K_G(\pt) = R(G)$.
Atiyah's definitions, applied to representations, are given as follows:
\begin{eqnarray*}
  P_n\negmedspace: R(G)&\to&R(G\times\Sn)\\
                   {[V]}&\mapsto& [V^{\tensor n}],
\end{eqnarray*}
is called the {\em $n^{th}$ power operation}\footnote{Atiyah's power
  operations actually take values in $R(G\wr\Sn)$ rather than
  $R(G\times\Sn)$, but for our 
  purposes, it is enough to pull them back via the diagonal map of $G^n$.}. 
For finite groups $G$ and $H$, and $h\in H$, Atiyah and Segal defined the map
\begin{eqnarray*}
  \tr h-\negmedspace : R(G\times H) &\to& R(G)\tensor_\ZZ\ZZ[\zeta]\\
                        W&\mapsto&\sum W_\eta\tensor\eta.      
\end{eqnarray*}
Here the sum runs over all eigenvalues $\eta$ of the action of
$h$ on $W$, the space $W_\eta$ is the eigenspace corresponding to $\eta$, and
$\zeta$ is a primitive $|h|^{th}$ root of unity.

The {\em $n^{th}$ symmetric power} of $V$ is given by
\begin{eqnarray*}
  \on{sym}_n(V) &=& \(P_n(V)\)^\Sn\\
                &=&\frac1{n!}\sum_{\sigma\in\Sn}\trs{P_n(V)}\\
                &=&\sum_{[\sigma]_\Sn}\GG{C_\sigma}\cdot\trs{P_n(V)}.
\end{eqnarray*}
Here $(-)^\Sn$ picks out the $\Sn$-invariant part of a
$G\times\Sn$-representation. All equalities are in $R(G)$, and the
second equality can be proved using the splitting principle.
Similarly, the {\em $n^{th}$ exterior power} of $V$ is given by picking out
the summand on which $\Sn$ acts by the sign representation:
\begin{eqnarray*}
  \lambda_n(V) &=& \Hom_{G\times\Sn}({\CC[G]}\tensor\sgn,P_n(V)).      \\
 &=&\frac1{n!}\sum_{\sigma\in\Sn}\sgn(\sigma)\cdot\trs{P_n(V)}.
\end{eqnarray*}

To define the $n^{th}$ Adams operation, let $c_n\in \Sn$ be a long
cycle, and set 
$$
  \psi_n(V) = \tr{c_n}{P_n(V)}.
$$
Again, one can use the splitting principle to compare this definition
with others in the literature.
The $\lambda_i$ make $R(G)$ into a $\Lambda$-ring. Writing 
$\Lambda_t=\sum\lambda_nt^n$ for the total exterior power, and
similarly $S_t$ for the total symmetric power, one has the well known
equations 
\begin{equation}
  \label{total-symmetric-power-1-Eqn}
  \Lambda_{-t}(x) =\exp\(-\sum\frac{\psi_n(x)}nt^n \)
  \quad\text{   and   }\quad
  S_t(x) =\exp\(\sum\frac{\psi_n(x)}nt^n\).   
\end{equation}
Here $t$ is a dummy variable. In \cite[9]{Ganter:thesis}, I explained
that, in the context of Atiyah's definitions, the equations
\eqref{total-symmetric-power-1-Eqn} are generating 
functions, which encode the cycle decomposition of the elements of the
symmetric groups. I also explained that for pairs of commuting
elements of the symmetric group, there is an analogous argument using
the decomposition of $\underline n$ into $\st$-orbits, i.e., into
transitive $\ZZ$-sets.

To make
the analogy with our setup even more explicit, we
will view class functions on $G$ as functions on $G/G$, where $G$ acts
on itself by conjugation, and note that $G/G$ is the coarse moduli
space of principal $G$-bundles on a circle. In this setup,
the $n^{th}$ Adams 
operation is given by pull-back along the degree $n$ self-map of the
circle. The group of deck transformations of such a map has order $n$.

The $n^{th}$ symmetric power is the weighted sum of pull-backs
along all isomorphism classes of 
$n$-fold coverings of $\bbS^1$. The weights are one over the order of
the automorphism groups.

From this point of view, the generating function equation for the
total symmetric power $S_t$ in 
\eqref{total-symmetric-power-1-Eqn} 
encodes the fact 
that every covering of the circle decomposes into the disjoint union
of covers with connected total space.
More precisely, 
the analogue of Proposition \ref{Sn-decomp-Prop} for covers
of the circle is
$$
  \coprod_{k\geq 0}\on{Cov}_k(\bbS^1) t^k = 
  \coprod_{m\geq 0}\(\coprod_{n\geq
    1}\on{Cov}_n^{\on{conn}}(\bbS^1)t^n\)\wr\Sigma_m.
$$
Here $\on{Cov}_k(\bbS^1)$ 
stands for
the groupoid of $k$-fold covers of $\bbS^1$ 
and $\on{Cov}^{\on{conn}}_k(\bbS^1)$ for the groupoid of those
with connected total
space.
The second equation in \eqref{total-symmetric-power-1-Eqn} is the weighted
sum of the points of both sides of 
this decomposition.
\subsection{Symmetric and exterior powers on $\Ell^\alpha_G$}
In \cite{Ganter:thesis}, I explained how Atiyah's theory can be
carried over to $\Ell_G$. The idea is simple: wherever Atiyah uses
elements of the group, one replaces them by pairs of commuting
elements. Then the Hecke operators take the role of the Adams
operations.
In this section, I will define the analogues of the symmetric and
exterior powers. We will denote them by $\sym{}$ and
$\lambda^{(2)}_n$. Depending on the taste of the reader, the $2$ can
stand for ``pairs of commuting elements'', for the fact that the torus is
made from two circles, or for ``chromatic level 2''. 

In order to define symmetric and exterior powers on $\Ell^\alpha_G$,
we need to 
extend the ``Hecke 
correspondence'' \eqref{Hecke-correspondence-Eqn} to all of $\mM_{G\times\Sn}$. 
We set
$$
  \map{t:= \mM_{id\times p}}{\mM_{G\times\Sn}}{\mM_G}.
$$
In order to define $s$, we recall the equivalence $e$ of Proposition
\ref{Sn-decomp-Prop}. For an abelian group $A$, let 
$$
  \map{s_A}{\mM_G\times_{\mM_1}\MzeroA}{\mM_G}
$$
be the source map.
We define $s$ to be $e$ composed with the map which on the component 
indexed by the decomposition $\sum N_A|A|=n$ is given by
$$
  \prod_A s_A\wr\Sigma_{N_A}.
$$
Over this component, we have
$$
  s^*\(\prod_A \mL^\alpha\wr\Sigma_{N_A}\)\cong
  t^*(\mL^\alpha)^{\tensor n}.
$$
\begin{Def}
  We define the $n^{th}$ symmetric power on sections $f$ of $\mL^{\alpha}$ by
  $$
    \sym n(f):= t_!\circ s^*
    \(\coprod_{n=\sum |A|N_A}\prod_A f\wr\Sigma_{N_A}\).
  $$
  Here $\sym n(f)$ is a section of $(\mL^\alpha)^{\tensor n}$.
  We define the total symetric power $\Sym q(f)$ by 
  $$\Sym q(f) := \sum_{n\geq 0}\sym n(f)\cdot q^n.$$
  Here $q$ is a dummy variable.
\end{Def}
From now on, $t$ will again denote our dummy variable (as opposed to
the target map).
The total symmetric power $\Sym t$ takes values in the graded ring
$$
  \bigoplus_{n\geq 0}\Ell_G^{n\alpha}(\pt)t^n.
$$
In the case that $\alpha=0$, this is the ring of formal power series
in $t$. The series $\Sym t(f)$ has a multiplicative inverse in 
$\bigoplus_{n\geq 0}\Ell_G^{-n\alpha}(\pt)t^n,$ which we will refer to
as the total 
(alternating) exterior power, denoted
$$
  \Lambda^{(2)}_{-t}(f) = \sum_{n\geq 0}\lambda_n^{(2)}(f)(-t)^n.
$$
The $\lambda^{(2)}_n$ make $\bigoplus\Ell^{-n\alpha}_G(\pt)$ a $\Lambda$-ring.
\begin{Def}
  Let $\st$ be a commuting pair in $\Sn$, and assume that $\st$
  induces the decomposition $n=\sum|A|N_A$.
  Let $s_\st$ denote the restriction of $s$ to the component of
  $\mM_\Sn$ corresponding to this decomposition. 
  We define the power operation $\psi_\st$ by
  $$
    \psi_\st(f) := s_\st^*(f).
  $$
\end{Def}
\begin{Def}
  Let $(\sigma,\rho)$ be a pair of commuting elements in $\Sn$. We
  define the signature $\sgn(\sigma,\rho)$ of the pair by
  $$\sgn(\sigma,\rho):=(-1)^{\# \text{even orbits of }(\sigma,\rho)}.$$
\end{Def}
We have
$$
  \sym n(f) = \frac1{n!} \sum_{\sigma\rho=\rho\sigma}
  \psi_{\st}(f)
$$
and 
$$
  \lambda_n^{(2)}(f) = \frac1{n!} \sum_{\sigma\rho=\rho\sigma}\sgnst
  \cdot \psi_{\st}(f). 
$$
The proof of \cite[Prop.9.1]{Ganter:thesis} goes through to give the
generating function equations
$$
  \Sym t(f) 
  = \exp \left[
    \sum_{k\geq 0} T_k(f)\cdot t^k
  \right]
$$
and
$$
  \Lambda^{(2)}_t(f)  
  = \exp \left[-
    \sum_{k\geq 0} T_k(f)\cdot t^k
  \right].
$$
\begin{Rem}
  Another proof for the first equation and a detailed discussion of the
  generating function argument in terms of covers over elliptic curves
  can be found in \cite{Roth}. In the
  case that $G$ is the trivial group and $f=1$, the statement is a special
  case of Roth's Lemma 2.9. His discussion is more general: it also applies
  to (simply branched) coverings of elliptic curves by curves of
  higher genus.
\end{Rem}
\subsection{Replicability}\label{replicability-Sec}
Using these definitions, we can rephrase Definition
\ref{untwisted-replicability-Def} as follows (compare 
\cite[(8.2)]{Borcherds}): 
\begin{Def}
  A McKay-Thompson series $f(q)$ is called {\em
    replicable}, if (over the untwisted sector) it satisfies
  \begin{equation}
    \label{twisted-replicability-Eqn}
    (f(t)-f(q)) = t\inv\Lambda_{-t}^{(2)}(f(q)).  
  \end{equation}
\end{Def}
A few words of explanation are in order. 
The left-hand side of (\ref{twisted-replicability-Eqn}) equals
$$
  t\inv-f(q)+\sum_{n\geq 1}a_nt^n,
$$
while the right-hand side is 
$$
  t\inv - f(q) + \sum_{n\geq 1}\lambda_{n+1}^{(2)} (-t)^n.
$$
Hence replicability means that for every $n\geq1$, the function
$\lambda^{(2)}_{n+1}(f)$ is a constant, namely $(-1)^{n+1}$ times the
$n^{th}$ coefficient of the $q$-expansion of $f$.
\section{The Witten genus}
Let $X$ be a compact, differentiable, spin manifold of
dimension $2d$, let $p_1$ denote the first Pontrjagin class, and
assume that $\frac{p_1}2(X)=0$.  
Let $T_\CC$ denote the complexification of 
the tangent bundle of $X$, and set $\overline T_\CC:=T_\CC-\CC^{2k}$. 
Then the Witten genus of $X$ is defined as
$$
  \Phi_\Witten(X) = \hat A\(X,\bigotimes_{k=1}^\infty
  \on{Sym}^{}_{q^k}(\ol T_\CC)\).
$$
Here $\on{Sym_q} = \on{Sym}_q^{(1)}$ stands for the total symmetric power in
$K$-theory, and $\hat A$ is the $\hat A$-genus. If the dimension of
$X$ is a multiple of $24$, the condition $\frac{p_1}{2}=0$ (often
called ``the vanishing of the {\em anomaly}'') implies
that  
\begin{equation}
  \label{A-hat-Eqn}
  \frac{\Phi_\Witten(X)}{\Delta^{\frac{2d}{24}}} =
  q^{-\frac{2d}{24}}\cdot \hat A\(X,\bigotimes_{k=1}^\infty
  \on{Sym}^{}_{q^k}(T_\CC)\)
\end{equation}
becomes a modular function. The Witten genus is linked to Moonshine
by Hirzebruch's prize question \cite[p.86]{Hirzebruch:Berger:Jung},
which was answered affirmatively by 
Hopkins and Mahowald \cite{Hopkins:Mahowald}:

\medskip
{\em Is there a $24$-dimensional manifold $X$ as above such that $\hat
  A(X)=1$ and $\hat A(X,T_\CC)=0$?}
\medskip

For such an $X$, the expression \eqref{A-hat-Eqn} becomes
$j-744$. Hirzebruch explains that if one could find an action of a finite
group $G$ on $X$ by diffeomorphisms, this
action would lift to the tangent bundle and its symmetric powers and
make \eqref{A-hat-Eqn} a McKay-Thompson series.
In the case that $G$ is the monster, one could hope for this series to
be that of classical Moonshine. 
In fact, assume that the action of $G$ lifts to the spin structure and that
$\(\frac{p_1}2(X)\)_G=0$. Then \eqref{A-hat-Eqn} is the $g=1$ part 
of a generalized McKay Thompson series, $\phi_{\Witten,G}(X)$, called the
{\em equivariant Witten genus} of $X$ \cite{Ganter:stringypowers}. 
Note that the equivariant genus in \cite{Ganter:stringypowers} differs
from the one usually found in the literature by a normalization
factor, so that my $\Phi_{\Witten}$ there is actually given by the
right-hand side of \eqref{A-hat-Eqn}.

However, the setup of \cite{Ganter:stringypowers} describes only the
$\alpha=0$
case. Also, the condition that the equivariant class 
$$\(\frac{p_1}2\)_G\in H^4(X\times_GEG;\ZZ)$$ 
should be zero turns out to be very restrictive. 
The following conjecture was explained to me by Matthew Ando, as was
much of the next section:
\begin{Conj}
  Let $X$ be a compact, differentiable, $G$-equivariant spin manifold
  of dimension $24d$, and assume that $\frac{p_1}2(X)=0$ and 
  $$\(\frac{p_1}2\)_G(X) = \pi^*(\alpha),$$ 
  where $\pi$ is the map from $X$ to a point, and $\alpha$ is an
  element of $H^4(BG;\ZZ)$. Then the equivariant Witten genus of $X$
  is a generalized Thompson series whose transformation behavior under
  $\SL$ makes
  it a section of $\mathcal L^{-\alpha}$.
\end{Conj}
The idea to interpret non-trivial anomaly as giving a section of a
line bundle goes back to Witten and can be found in the work of 
Freed and Witten \cite{Freed:Witten},
Liu \cite{Liu}, and Ando \cite{Ando:circle}.
The next section outlines how this conjecture could follow
from the conjectural behavior of the Thom-isomorphism in
\cite{Ginzburg:Kapranov:Vasserot}. 
\subsection{Twisted Thom isomorphisms}
The Witten genus is closely related to the theory of Thom isomorphisms
and transfer maps in elliptic cohomology. Assume that $V-TX$ is a spin
bundle and that we are given a lift of the action of $G$ to the spin
structure. Let $X^{V-TX}$ denote the equivariant Thom spectrum of the
virtual vector bundle $V-TX$ \cite{Lewis:May:Steinberger}.
Then $\Ell^0_G(X^{V-TX})$ is an invertible $\Ell^0_G$ module
sheaf. Pulling back first along the relative zero section $X^{-TX}\to
X^{V-TX}$ and then along the Pontrjagin-Thom collapse $\bbS^0\to X^{-TX}$, 
we obtain a map
$$
  \Ell^0_G(X^{V-TX})\to  \Ell^0_G(X^{-TX}) \to
  \Ell^0_G(\bbS^0)=\mathscr O_{\mM_G}.
$$
If the equivariant characteristic class 
$$\(\frac{p_1}{2}\)_G(V-TX)\in H^4(EG\times_GX;\ZZ)$$
equals zero, the Thom isomorphism is a trivialization
$$
    \Ell^0_G(X)\stackrel\cong\longrightarrow  \Ell^0_G(X^{V-TX}).
$$
Composing it with the map above, gives the {\em transfer-} (or {\em
  Gysin-}) map 
along $\map\pi {X_+}{\bbS^0}$,
$$
  \map{\pi_!^V}{\Ell^0_G(X)}{\Ell^0_G(\bbS^0)}, 
$$
and the Witten genus of $X$ twisted with $V$ is
$$
  \phi_{\Witten,G}(V) = \pi_!^V(1). 
$$
If $V$ is trivial, we drop it from the notation and have
$$\phi_{\Witten,G}(X) = \pi_!(1).$$ 
Assume now that the weaker condition
$$\(\frac{p_1}{2}\)_G(V-TX) = \pi^*(\alpha)$$
is satisfied, where
$$\alpha\in H^4(BG;\ZZ),$$
and assume that, non-equivariantly, we have 
$$\(\frac{p_1}{2}\)(V-TX) = 0\in H^4(X;\ZZ).$$
Then we get a twisted Thom isomorphism
$$
  \Ell^0_G(X)\tensor_{\mathscr O_{\mM_G}}\mL^{\alpha}\stackrel\cong\longrightarrow
  \Ell^0_G(X^{V-TX}), 
$$
yielding a transfer map
$$
  \map{\pi_!^V}{\Ell^0_G(X)}{\(\mL^\alpha\)\inv}.
$$
Hence the Witten genus $\phi_{\Witten,G}(V)$ takes values in
the global sections of $\mL^{-\alpha}$. 
\subsection{A few words about physics}
Principal bundles over Riemann surfaces turn up in string theory, when the
target space of the theory is an orbifold quotient $X\mmod G$. 
In string theory with target space $X$, one considers spaces of maps
from the circle (a {\em closed string}) or from a Riemann surface $\Sigma$
(the {\em worldsheet}) to $X$.
An orbifold map from $\Sigma$ to $X\mmod G$ turns out to be the same as an
equivalence class of $G$-equivariant maps from a principal $G$-bundle
over $\Sigma$ to $X$: 
$$
  \on{map}_\orb(\Sigma, X\mmod G)= 
  \coprod_{\stackrel{P\to \Sigma}{\text{pbdl}}}\on{map}_G(P,X)\slash \sim 
$$
(cf., e.g., \cite{Sharpe:Quotient}, \cite{Moerdijk}, \cite{Tamanoi},
\cite{Lupercio:Uribe:Xicotencatl}).  
The space 
$$\mL(X\mmod G) = \on{map}_\orb(\bbS^1,X\mmod G)$$
is called the {\em orbifold loop space} of $X\mmod G$, and the space
$$\mL^2(X\mmod G)=\on{map}_\orb(\TT^2,X\mmod G)$$ is the double (orbifold) loop
space of $X\mmod G$. They decompose into components 
$$
  \mL_g(X\mmod G) := \on{map}_G(P_g,X) \quad\text{and}\quad
  \mL^2_{g,h}(X\mmod G) = \on{map}_G(P_{g,h},X),
$$
often referred to as the {\em twisted sectors} and the {\em membrane twisted
sectors} of the theory.

In the non-equivariant case, the Witten genus $\Phi_\Witten(X)$ is the
genus $1$ partition function of the supersymmetric non-linear sigma model for
the target space $X$ \cite{Stolz:Teichner}.

In \cite{deFernex:Lupercio:Nevins:Uribe}, de Fernex, Lupercio, Nevins
and Uribe consider
a supersymmetric string sigma
model whose 
target space is an orbifold and show that its partition function
on a two-dimensional torus
is the orbifold Witten genus. The orbifold Witten genus is
related to the equivariant Witten genus by
\begin{eqnarray*}
  \Phi_{\Witten,\orb}(X\mmod G) &=& \GG G\sum_{gh=hg}\Phi_{\Witten, G}(X)(g,h)\\
  &=&   \ind G1\(\Phi_{\Witten, G}(X)\).
\end{eqnarray*}
In \cite{deFernex:Lupercio:Nevins:Uribe}, the $(g,h)$-summand of
$\Phi_{\Witten,\orb}$ is the contribution of the corresponding twisted
sector to the partition sum,  
which is a phase integral over the torus 
with boundary conditions twisted by $g$ and $h$.
Hence the discussion of \cite{deFernex:Lupercio:Nevins:Uribe} fits the
equivariant Witten genus (the same thing as in
\cite{Ganter:stringypowers}) into the framework of non-linear sigma models.


The idea to study maps into an orbifold target space via
principal bundles on the worldsheet already
turned up in Segal's celebrated work
\cite{Segal:CFT}. In the same paper, Segal prominently announced the
existence of a ``{\em fairly simple and natural conformal field
  theory}'', whose 
automorphism group is the Monster.
In \cite{Frenkel:Lepowsky:Meurman}, 
Frenkel, Lepowski and Meurman constructed a vertex algebra
whose automorphism group is the Monster, the so-called {\em Moonshine module}.
Vertex algebras are often thought of as the mathematical language of
two-dimensional conformal field theory.
In his famous paper \cite{Borcherds}, Borcherds proved that the
Moonshine module is indeed the McKay-Thompson series of the classical
Moonshine conjecture. It is expexted and in many cases proved
that generalized Moonshine is given by the twisted sectors of
the Moonshine module
(cf.\ \cite{Dong:Li:Mason:modinv}).  

In \cite{Beautyandthebeast}, Dixon, Ginsparg, and Harvey translate the
result of \cite{Frenkel:Lepowsky:Meurman} into the language of string theory.
According to them,
the target space of the Monster CFT is the orbifold 
$$\TT^{24}_{\on{Leech}}\mmod (\ZZ/2\ZZ).$$
Here $\TT^{24}_{\on{Leech}} = \mathbb R^{24}/\Lambda_{\on{Leech}}$,
where $\Lambda_{\on{Leech}}$ is the Leech lattice, and the $\ZZ/2\ZZ$
action is 
induced by multiplication with $-1$ on $\mathbb R^{24}$.
Its orbifold partition function is $j-744$. The
automorphism group of the theory is the
Monster. For a pair of commuting elements $(g,h)$ in the Monster, 
\cite{Beautyandthebeast} describes the generalized Moonshine
function $f(g,h;\tau)$ as the phase 
integral over the torus with boundary conditions twisted by $g$ and
$h$. 

Given these facts, one could hope to reinterpret Hirzebruch's question
and ask whether the 
orbifold Witten genus of $\TT^{24}_{\on{Leech}}\mmod(\ZZ/2)$ was given by
$j-744$. After all, the Witten genus is the partition function of the
theory, and the partition function is $j-744$. Unfortunately, the two
theories are different, and in particular, their 
notions of (orbifold) partition function 
do not agree. For instance, as a consequence of the rigidity theorem,
the contribution of the non-twisted sector of the Witten genus,
$\Phi_\Witten(\TT_{\on{Leech}}^{24})$, is zero, whereas the analogous
contribution to the partition function of the Moonshine module is $j-720$.

When Witten defines the Witten genus in \cite{Witten:ellipticgenera},
he explains that different non-linear sigma models lead to different genera.
\begin{Que}
  Is it possible to describe the partition function 
  of the ``holomorphic conformal field theory'' in
    \cite{Beautyandthebeast} as a push-forward construction in
    equivariant elliptic cohomology, similar to the description of the
    Witten genus in \cite{Ganter:stringypowers}, and can we identify the
    corresponding genus? 
    How does the discussion in \cite{deFernex:Lupercio:Nevins:Uribe}
    need to be modified to supply the geometric side of this
    construction (in terms of orbifold loop spaces)? 
    Ideally, one would want this genus to preserve power operations. In
    this case, the theory would have a
    Dijkgraaf-Moore-Verlinde-Verlinde formula.
\end{Que}
\subsection{The Witten genus and replicability?}
In \cite{Ganter:stringypowers}, I calculated the total symmetric power in Devoto's
equivariant Tate $K$-theory applied to a $G$-equivariant vector bundle
$V$ and observed that it equals Witten's exponential characteristic class:
$$
  \on{Sym}_q^{(2)}(V) = \bigotimes_{k=1}^\infty \on{Sym}_{q^k}^{(1)}(V).
$$
Recall that we have
  $$\on{Sym}_q^{(2)}(V) = \Lambda^{(2)}_{-q}(-V).$$
\begin{Que}
  Can one describe replicability
  \eqref{twisted-replicability-Eqn} in the context of Hirzebruch's prize
  question as asking for an equality about iterated total exterior
  powers at chromatic level $2$?
\end{Que}
\begin{Que}
  Is it possible to make sence of the notion of replicability 
  for generalized McKay Thompson series? In this case,
  the coefficient $\lambda^{(2)}_{n+1}(f(q))$ 
  on the right-hand side of \eqref{twisted-replicability-Eqn}, 
  is a section of the line bundle $(\mL^{\alpha})^{\tensor n+1}$. What does it
  mean for such a thing to be a constant? And how should we take the
  fractional powers of $t$ occurring over the twisted sectors on the
  left-hand side   
  into account? 
\end{Que}
Note that these two issues are related:
  The fractional powers in the $q$-expansions of the twisted sectors
  come from the fact that $T$ does not act 
  trivially on the lines $L_{g,h}$.
  Let $L_N$ be the line of invariant sections of the functor $\mathcal
F$ on $\mu_N\rtimes\mu_N$ which is given by
$
 \mathcal F(\zeta) = \CC, 
$
on objects and such that 
$$\map{\mathcal F(\zeta,\xi)}{\mathcal F(\zeta)}{\mathcal
  F(\zeta\xi)}$$ 
is multiplication with $\xi$.
Then there is a $\ZZ/N\ZZ$-equivariant ring isomorphism
\begin{eqnarray*}
  \bigoplus_{m\geq -1} L_N^{\tensor m}q^m&\to&\CC\ps{q^\frac1N}  \\
  a_mq^m&\mapsto& a_m(e^{\frac{2\pi i}N})q^\frac mN.
\end{eqnarray*}
Here $1\in\ZZ/N\ZZ$ acts on $q^\frac1N$ by multiplication with
$e^{\frac{2\pi i}N}$.

%
%
\appendix
\section{Generalities on principal bundles}
This appendix collects some well known facts about principal bundles
used in the paper.
Let $G$ be a finite group.
\begin{Def}
  A principal $G$-bundle $\xi$ consists of a map $\map \pi PX$, where
  $P$ is a right $G$-space such that $G$ acts strictly transitively on
  the fibers of $\pi$. The space $P$ is called the {\em total space} and $X$ 
  the {\em base space} of $\xi$.
\end{Def}
For simplicity, we will always assume that the base space $X$
is connected.
\begin{Def}
  Let $i\negmedspace :G\hookrightarrow H$ be an inclusion of
  groups, and let $\xi$ be a principal $G$-bundle.
  The principal $H$-bundle $\xi[H]$ associated to $\xi$ via $i$ has
  total space 
  $$P\times_GH:=(P\times H)/\sim,$$
  where $(p,h)\sim (pg,i(g^{-1})h)$. The action of $H$ on $P\times_GH$
  is defined by $$(p,h)\cdot h':=(p,hh').$$ The projection $\pi[H]$ is
  defined by $\pi[H](p,h):=\pi(p)$.
\end{Def}
If $\map fX\BG$ classifies the $G$-bundle $\xi$,
the associated $H$-bundle $\xi[H]$ is classified by the map $\on{Bi}\circ f$.
\begin{Def}
Let $\xi=(\map \pi PX)$ be a principal $G$-bundle with connected total space,
and let $\map\alpha HG$ be an isomorphism of groups. We define the
pull-back of $\xi$ along $\alpha$ to be the principal
$H$-bundle $\alpha^*(\xi)$ consisting of the same map $\pi$, where the
$H$ action on the total space $P$ is given by
$$p*h:=p\cdot\alpha(h).$$ In the case $H=G$, we will refer to this
pulled back action as ``{\em the $G$-action
on $P$ twisted by $\alpha$}''.
\end{Def}
\begin{Lem}
  The bundle $\alpha^*(\xi)$ is isomorphic to the $G$-bundle
  associated to $\xi$ via the inclusion $\map{\alpha^{-1}}GG$.
\end{Lem}
\begin{Pf}{}
  Consider the map
  \begin{eqnarray*}
    f\negmedspace : P&\to&P\times_{\alpha^{-1}}G\\
    p&\mapsto&(p,1).
  \end{eqnarray*}
  It is obviously fiber preserving and bijective.
  We have
  $$f(p)\cdot g = (p,1)\cdot g = (p,g) \sim
  (p\alpha(g),\alpha^{-1}(\alpha(g))\cdot 1) =
  (p\alpha(g),1)=f(p\alpha(g))=f(p*g).$$
  Therefore $f$ is $G$-equivariant with respect to the $G$-action on
  $P$ twisted by $\alpha$.
\end{Pf}
\begin{Cor}\label{Balpha-Cor}
  Let $\xi$ be classified by the map $\map fX\BG$. Then
  $\alpha^*(\xi)$ is classified by
  $$\on{B\alpha^{-1}}\circ f.$$
\end{Cor}
\begin{Lem}\label{associated-fiber-Lem}
  Let $G\sub H$ be an inclusion of groups, let $\xi$ be a principal
  $H$-bundle,  
  and let $F$ be a left
  $G$-space.
  Then
  $$\xi\times_HF=\xi[G]\times_GF.$$
\end{Lem}
%
%
%
\begin{Prop}\label{GHisom-Prop}
  Let $\xi=(\map \pi PX)$ be a principal $G$-bundle over $X$, and
  $\zeta=(\map \rho QX)$ be a principal $H$-bundle over $X$, and assume
  that both total spaces $P$ and $Q$ are connected. Assume further 
  that there exists an isomorphism of covering spaces of $X$
  $$\map fPQ.$$
  Then there exists an isomorphism of groups
  $$\map i GH$$
  making $f$ into an isomorphism of principal $G$-bundles
  $$\isomap f\xi{ i^*\zeta}.$$
  Moreover, $ i$ is uniquely determined by $f$. 
\end{Prop}
\begin{Pf}{}
  Fix an element $p\in P$.
  We define the map $\map i GH$ as follows: for $g\in G$ let
  $ i(g)$ be the unique element of $H$ mapping $f(p)$ to $f(pg)$.
  It is clear that $ i$ is a bijection. We need to show that
  $ i$ is a map of groups. Fix $g\in G$, and consider the
  deck transformations $D_g$ of $\pi$ sending $p$ to $pg$ and
  $D_{ i(g)}$ of $\rho$ sending $q$ to  $q\cdot i(g)$.
  Then $D_{ i(g)}$ and $f\circ D_g\circ f^{-1}$ both are
  deck transformations of $Q$ sending $f(p)$ to $f(pg)$.
  We may conclude that they are the same deck transformation. 
  Now we have
  \begin{eqnarray*}
    f(p)\cdot i(hg) &=& f(p(hg))\\
    &=& f((f^{-1}(f(ph)))\cdot g) \\
    &=& f(ph)\cdot i(g)\\
    &=& (f(p)\cdot i(h))\cdot i(g).
  \end{eqnarray*}
  Further, for every $g\in G$, we have
  $$f\circ D_g=D_ i(g)\circ f,$$
  proving that $f$ is an isomorphism of principal $G$-bundles.
  Let now $\map\beta GH$ be an isomorphism making $f$ into an
  isomorphism between the principal $G$-bundles $\xi$ and
  $\beta^*(\zeta)$.
  Then we have for every $g\in G$,
  $$f(pg)=f(p)\beta(g),$$
  implying that $\beta(g)= i(g)$.
\end{Pf}
%
%

\bibliographystyle{alpha}
\bibliography{twistedhecke.bib}
\end{document}